\documentclass[onecolumn,twoside,10pt]{IEEEtran}
\usepackage{amsmath}
\usepackage{amssymb}
\usepackage{url}
\usepackage{mathtools}
\usepackage{color}
\usepackage{algorithm}
\usepackage[noend]{algorithmic}
\usepackage{graphicx}
\usepackage{subfigure}
\usepackage{tabularx,booktabs}

\newtheorem{thm}{Theorem}
\newtheorem{rem}[thm]{Remark}
\newtheorem{prop}[thm]{Proposition}

\newcommand{\R}{\mathbb{R}}

\newcommand{\C}{\mathbb{C}}
\newcommand{\cX}{\mathcal{X}}
\newcommand{\cD}{\mathcal{D}}
\newcommand{\cS}{\mathcal{S}}
\newcommand{\cF}{\mathcal{F}}
\newcommand{\cN}{\mathcal{N}}
\newcommand{\cE}{\mathcal{E}}
\newcommand{\cR}{\mathcal{R}}

\renewcommand{\l}{\ell}
\newcommand{\proj}{\mathcal{P}}
\newcommand{\suchthat}{\text{s.t.}}
\newcommand{\norm}[1]{\lVert {#1} \rVert}
\newcommand{\bignorm}[1]{\big\lVert {#1} \big\rVert}
\newcommand{\Norm}[1]{\left\lVert {#1} \right\rVert}
\newcommand{\abs}[1]{\lvert {#1} \rvert}

\newcommand{\st}{~~\text{s.t.}~~}
\newcommand{\define}{\coloneqq}
\newcommand{\NP}{\textsf{NP}}
\renewcommand{\phi}{\varphi}

\DeclareMathOperator{\sign}{sign}

\DeclareMathOperator{\vect}{vec}
\DeclareMathOperator{\midpt}{mid}

\algsetup{indent=2em}

\newcommand{\F}{\mathbf{F}}
\newcommand{\M}{\mathbf{M}}

\newcommand{\A}{\mathbf{A}}
\renewcommand{\P}{\mathbf{P}}
\newcommand{\B}{\mathbf{B}}
\newcommand{\X}{\mathbf{X}}
\newcommand{\Y}{\mathbf{Y}}
\newcommand{\G}{\mathbf{G}}
\renewcommand{\H}{\mathbf{H}}
\newcommand{\D}{\mathbf{D}}
\newcommand{\x}{\mathbf{x}}
\renewcommand{\d}{\mathbf{d}}
\newcommand{\y}{\mathbf{y}}

\newcommand{\z}{\mathbf{z}}
\renewcommand{\a}{\mathbf{a}}
\newcommand{\al}{\mathbf{\alpha}}
\newcommand{\I}{\mathbf{I}}

\newcommand{\mZ}{\mathbf{Z}}
\newcommand{\mN}{\mathbf{N}}
\newcommand{\mC}{\mathbf{C}}
\newcommand{\n}{\mathbf{n}}
\newcommand{\mR}{\mathbf{R}}
\newcommand{\q}{\mathbf{q}}

\newcommand{\zeros}{\mathbf{0}}

\begin{document}

\title{DOLPHIn -- Dictionary Learning for\\Phase Retrieval}%


\author{Andreas M. Tillmann,~\IEEEmembership{} Yonina
  C. Eldar,~\IEEEmembership{Fellow,~IEEE,}
  and~Julien~Mairal,~\IEEEmembership{Member,~IEEE}
  \thanks{A. M. Tillmann is with the TU Darmstadt, Research Group
    Optimization, Dolivostr. 15, 64293 Darmstadt, Germany (e-mail:
    tillmann@mathematik.tu-darmstadt.de).}%
  \thanks{Y. C. Eldar is with the Department of Electrical
    Engineering, Technion -- Israel Institute of Technology, Haifa
    32000, Israel (e-mail: yonina@ee.technion.ac.il).}%
  \thanks{J. Mairal is with Inria, Lear Team, Laboratoire Jean
    Kuntzmann, CNRS, Universit\'e Grenoble Alpes, 655, Avenue de
    l'Europe, 38330 Montbonnot, France (e-mail:
    julien.mairal@inria.fr).}  
	\thanks{The work of J. Mairal was funded by the French 
	National Research Agency [Macaron project, ANR-14-CE23-0003-01]. The
	work of Y. Eldar was funded by the European Union's Horizon 2020
	research and innovation programme under grant agreement ERC-BNYQ, and
	by the Israel Science Foundation under Grant no. 335/14.}%
  \thanks{This work has been submitted to
    the IEEE for possible publication. Copyright may be transferred
    without notice, after which this version may no longer be
    accessible.}}


\maketitle

\begin{abstract}
  We propose a new algorithm to learn a dictionary for reconstructing
  and sparsely encoding signals from measurements without
  phase. Specifically, we consider the task of estimating a
  two-dimensional image from squared-magnitude measurements of a
  complex-valued linear transformation of the original image. Several
  recent phase retrieval algorithms exploit underlying sparsity of the
  unknown signal in order to improve recovery performance. In this
  work, we consider such a sparse signal prior in the context of phase
  retrieval, when the sparsifying dictionary is not known in
  advance. Our algorithm jointly reconstructs the unknown
  signal---possibly corrupted by noise---and learns a dictionary such
  that each patch of the estimated image can be sparsely
  represented. Numerical experiments demonstrate that our approach can
  obtain significantly better reconstructions for phase retrieval
  problems with noise than methods that cannot exploit such ``hidden''
  sparsity. Moreover, on the theoretical side, we provide a
  convergence result for our method.
\end{abstract}

\begin{IEEEkeywords}
  (
  MLR-DICT, MLR-LEAR, OPT-NCVX, OPT-SOPT)
  Machine Learning, Signal Reconstruction, Image Reconstruction
\end{IEEEkeywords}

\IEEEpeerreviewmaketitle

\section{Introduction}\label{sec:intro}
\IEEEPARstart{P}{hase} retrieval has been an active research topic for
decades~\cite{M07,SECCMS15}. The underlying goal is to estimate an
unknown signal from the modulus of a complex-valued linear
transformation of the signal.  With such nonlinear measurements, the
phase information is lost (hence the name ``phase retrieval''),
rendering the recovery task ill-posed and, perhaps not surprisingly,
\NP-hard~\cite{SC91}. Traditional approaches consider cases where the
solution is unique up to a global phase shift, which can never be
uniquely resolved, and devise signal reconstruction algorithms for
such settings. Uniqueness properties and the empirical success of
recovery algorithms usually hinge on oversampling the signal, i.e.,
taking more measurements than the number of signal components.

The most popular techniques for phase retrieval are based on
alternating projections, see \cite{GS72,F82,BCL02} for
overviews. These methods usually require precise prior information
about the signal (such as knowledge of the support set) and often
converge to erroneous results. More recent approaches include
semidefinite programming relaxations
\cite{CESV13,SESS11,JOH12,OYDS12,WdAM15} and gradient-based methods
such as Wirtinger Flow \cite{CLS14,CLM15}.

In recent years, new phase retrieval techniques were developed for
recovering \emph{sparse} signals, which are linear combinations of
only a few atoms from a known dictionary
\cite{SESS11,SBE14,CLM15,MRB07}. With a sparsity assumption, these
algorithms obtained better recovery performance than traditional
non-sparse approaches.  The main idea is akin to compressed sensing,
where one works with fewer (linear) measurements than signal
components \cite{FR13,EK12,E15}. An important motivation for
developing sparse recovery techniques was that many classes of signals
admit a sparse approximation in some basis or overcomplete
dictionary~\cite{OF96,EA06,MBP14}. While sometimes such dictionaries
are known explicitly, better results have been achieved by adapting
the dictionary to the data, e.g., for image denoising \cite{EA06}.
Numerous algorithms have been developed for this task, see,
e.g.,~\cite{OF96,AEB06,MBPS10}. In this traditional setting, the
signal measurements are~\emph{linear} and a large database of training
signals is used to train the dictionary.

In this work, we propose a dictionary learning formulation for
simultaneously solving the signal reconstruction and sparse
representation problems given \emph{nonlinear, phaseless} and
\emph{noisy} measurements. To optimize the resulting (nonconvex)
objective function, our algorithm---referred to as DOLPHIn (DictiOnary
Learning for PHase retrIeval)---alternates between several
minimization steps, thus monotonically reducing the value of the
objective until a stationary point is found (if step sizes are chosen
appropriately). Specifially, we iterate between best fitting the data and
sparsely representing the recovered signal. DOLPHIn combines projected
gradient descent steps to update the signal, iterative shrinkage to
obtain a sparse approximation \cite{BT09}, and block-coordinate
descent for the dictionary update \cite{MBPS10}.

In various experiments on image reconstruction problems, we
demonstrate the ability of DOLPHIn to achieve significantly improved
results when the oversampling ratio is low and the noise level high,
compared to the recent state-of-the-art Wirtinger Flow (WF) method
\cite{CLS14}, which cannot exploit sparsity if the dictionary is
unknown. In this two-dimensional setting, we break an image down into
small patches and train a dictionary such that each patch can be
sparsely represented using this dictionary. The patch size as well as
the amount of overlap between patches can be freely chosen, which
allows us to control the trade-off between the amount of computation
required to reconstruct the signal and the quality of the result.

The paper is organized as follows: In Sections~\ref{sec:probstatement}
and~\ref{sec:algo}, we introduce the DOLPHIn framework and
algorithm. Then, in Section~\ref{sec:experiments}, we present
numerical experiments and implementation details, along with
discussions about (hyper-)parameter selection and variants of DOLPHIn.
We conclude the paper in Section~\ref{sec:conclusion}. 
The appendices provide 
further details on the mathematical derivation of the DOLPHIn
algorithm and its properties. 
A short preliminary version of this work
appeared in the conference paper~\cite{TEM15}.

\section{Phase Retrieval Meets Dictionary Learning}\label{sec:probstatement}
In mathematical terms, the phase retrieval problem can be formulated
as solving a nonlinear system of equations:
\begin{equation}\label{eq:PR}
  \text{Find }\x\in\cX\subseteq\C^{N}\quad\text{s.t.}\quad \abs{f_i(\x)}^2=y_i~~~\forall\,i=1,\dots,M,
\end{equation}
where the functions $f_i:\C^N\to\C$ are linear operators and the
scalars $y_i$ are nonlinear measurements of the unknown original
signal $\hat\x$ in $\cX$, obtained by removing the phase
information. The set~$\cX$ represents constraints corresponding to
additional prior information about~$\hat\x$. For instance, when
dealing with real-valued bounded signals, this may typically be a box
constraint~$\cX=[0,1]^N$.  Other common constraints include
information about the support set---that is, the set of nonzero
coefficients of~$\hat\x$.  Classical phase retrieval concerns the
recovery of~$\hat\x$ given the (squared) modulus of the signal's
Fourier transform. Other commonly considered cases pertain to
randomized measurements ($f_i$ are random linear functions) or coded
diffraction patterns, i.e., concatenations of random signal masks and
Fourier transforms (see, e.g., \cite{SECCMS15,CLS14}).

\subsection{Prior and Related Work}
The most popular methods for classical phase retrieval---Fienup's
algorithm~\cite{F82} and many related approaches \cite{SECCMS15,GS72,BCL02,BCL03,L05}---are based on
alternating projections onto the sets
$\mathcal{Y}\define\{\x \in \C^N \st \abs{f_i(\x)}=y_i~\forall\,i\}$ (or
$\{\x \in \C^N \st \abs{f_i(\x)}^2=y_i~\forall\,i\}$) and onto
the set~$\cX$. However, the nonconvexity of $\mathcal Y$ makes the projection
not uniquely defined and possibly hard to compute. The
success of such projection-based methods hinges critically on precise
prior knowledge (which, in general, will not be available in practice)
and on the choice of a projection operator onto $\mathcal{Y}$. Ultimately,
convergence to $\hat\x$ (up to global phase) is in general not guaranteed
and these methods often fail in practice.

Further algorithmic techniques to tackle~\eqref{eq:PR} include two
different semidefinite relaxation approaches, PhaseLift~\cite{CESV13}
and PhaseCut~\cite{WdAM15}. PhaseLift ``lifts'' \eqref{eq:PR} into the
space of (complex) positive semidefinite rank-$1$ matrices via the
variable transformation $\X\define\x\x^*$. Then, the nonlinear
constraints $\abs{f_i(\x)}^2=y_i$ are equivalent to
\emph{linear} constraints with respect to  the matrix variable $\X$. By
suitably relaxing the immediate but intractable rank-minimization
objective, one obtains a convex semidefinite program (SDP). Similarly,
PhaseCut introduces a separate variable $\mathbf{u}$ for the phase,
allowing to eliminate~$\x$, and then lifts $\mathbf{u}$ to obtain an
equivalent problem with a rank-$1$-constraint, which can be dropped to
obtain a different SDP relaxation of~\eqref{eq:PR}. Despite some
guarantees on when these relaxations are tight, i.e., allow for
correctly recovering the solution to~\eqref{eq:PR} (again up to a
global phase factor), their practical applicability is limited due to
the dimension of the SDP that grows quadratically with the problem dimension.

A recent method that works in the original variable space is the
so-called \emph{Wirtinger Flow} algorithm~\cite{CLS14}. Here, \eqref{eq:PR}
is recast as the optimization problem
\begin{equation}\label{eq:WFobjective}
\min_{\x\in\C^N}\tfrac{1}{4M}\sum_{i=1}^M \big(\abs{f_i(\x)}^2-y_i\big)^2,
\end{equation}
which is approximately solved by a gradient descent algorithm.  Note
that in the case of complex variables, the concept of a gradient is
not well-defined, but as shown in~\cite{CLS14}, a strongly related
expression termed the ``Wirtinger derivative'' can be used instead and
indeed reduces to the actual gradient in the real case. For the case
of i.i.d. Gaussian random measurements, local convergence with high
probability can be proven for the method, and a certain spectral
initialization provides sufficiently accurate signal estimates for
these results to be applicable. Further variants of the Wirtinger Flow
(WF) method that have been investigated are the Truncated
WF~\cite{CC15}, which involves improving search directions by a
statistically motivated technique to filter out components that bear
``too much'' influence, and Thresholded WF~\cite{CLM15}, which allows
for improved reconstruction of \emph{sparse} signals (i.e., ones with
only a few significant components or nonzero elements), in particular
when the measurements are corrupted by noise.

The concept of sparsity has been successfully employed in the context
of signal reconstruction from \emph{linear} measurements, perhaps most
prominently in the field of \emph{compressed
  sensing}~\cite{FR13,EK12,E15,E10} during the past decade. There, the
task is to recover an unkown signal $\hat\x\in\C^N$ from $M<N$ linear
measurements---that is, finding the desired solution among the
infinitely many solutions of an underdetermined system of linear
equations. For signals that are (exactly or approximately) sparse with
respect to some basis or dictionary, i.e., when
$\hat\x\approx\D\hat\a$ for a matrix $\D$ and a vector $\hat\a$ that
has few nonzero entries, such recovery problems have been shown to be
solvable in a very broad variety of settings and applications, and
with a host of different algorithms. Dictionaries enabling sparse
signal representations are sometimes, but not generally, known in
advance. The goal of \emph{dictionary learning} is to improve upon the
sparsity achievable with a given (analytical) dictionary, or to find a
suitable dictionary in the first place.  Given a set of training
signals, the task consists of finding a dictionary such that every
training signal can be well-approximated by linear combinations of
just a few atoms. Again, many methods have been developed for this
purpose (see, e.g.,~\cite{OF96,EA06,MBP14,AEB06,MBPS10}) and
demonstrated to work well in different practical applications.

Signal sparsity (or compressability) can also be beneficially
exploited in phase retrieval methods,
cf.~\cite{SESS11,JOH12,CLM15,SBE14,MRB07}. However, to the best of our
knowledge, existing methods assume that the signal is sparse itself or
sparse with respect to a fixed pre-defined dictionary. This motivates
the development of new algorithms and formulations to \emph{jointly}
learn suitable dictionaries and reconstruct input signals from
nonlinear measurements.

\subsection{Dictionary Learning for Phase Retrieval}
In this paper, we consider the problem of phase retrieval by focusing on image
reconstruction applications. Therefore, we will work in a two-dimensional
setting directly. However, it should be noted that all expressions and
algorithms can also easily be formulated for one-dimensional signals
like~\eqref{eq:PR}, as detailed in 
Appendix~\ref{sec:appendixA}. 
We will also consider the case of noisy measurements, and will show that
our approach based on dictionary learning is particularly robust to noise,
which is an important feature in practice.

Concretely, we wish to recover an image $\hat\X$ in $[0,1]^{N_1\times N_2}$ from
noise-corrupted phaseless nonlinear measurements
\begin{equation}
   \Y\define\abs{\cF(\hat\X)}^2+\mN, \label{eq:prb}
\end{equation}
where $\cF:\C^{N_1\times N_2}\to\C^{M_1\times M_2}$ is a linear
operator, $\mN$ is a real matrix whose entries represent noise, and
the complex modulus and squares are taken component-wise. As mentioned
earlier, signal sparsity is known to improve the performance of phase
retrieval algorithms, but a sparsifying transform is not always known
in advance, or a better choice than a predefined selection can
sometimes be obtained by adapting the dictionary to the data. In the
context of image reconstruction, this motivates learning a dictionary
$\D$ in $\R^{s\times n}$ such that each $s_1\times s_2$ patch
$\hat\x^i$ of $\hat\X$, represented as a vector of size $s=s_1 s_2$,
can be approximated by $\hat\x^i\approx\D\a^i$ with a sparse vector
$\a^i$ in $\R^n$. Here, $n$ is chosen a priori and the number of
patches depends on whether the patches are overlapping or not. In
general, $\D$ is chosen such that $n\geq s$. With linear measurements,
the paradigm would be similar to the successful image denoising
technique of~\cite{EA06}, but the problem~(\ref{eq:prb}) is
significantly more difficult to solve due to the modulus operator.

Before detailing our algorithm for solving~(\ref{eq:prb}), we
introduce the following notation. Because our approach is patch-based
(as most dictionary learning formulations), we consider the linear
operator $\cE:\C^{N_1\times N_2}\to\C^{s\times p}$ that extracts
the~$p$ patches $\x^i$ (which may overlap or not) from an image $\X$
and forms the matrix $\cE(\X)=(\x^1,\dots,\x^p)$.  Similarly, we
define the linear operator $\cR:\C^{s\times p}\to\C^{N_1\times N_2}$
that reverses this process, i.e., builds an image from a matrix
containing vectorized patches as its columns. Thus, in particular, we
have $\cR(\cE(\X))=\X$.  When the patches do not overlap, the
operator~$\cR$ simply places every small patch at its appropriate
location in a larger $N_1 \times N_2$ image.  When they overlap, the
operator $\cR$ averages the pixel values from the patches that fall
into the same location.  Further, let $\A\define(\a^1,\dots,\a^p)$ in
$\R^{n\times p}$ be the matrix containing the patch representation
coefficient vectors as columns. Then, our desired sparse-approximation
relation ``$\x^i\approx\D\a^i$ for all $i$'' can be expressed as
$\cE(\X)\approx\D\A$.

With this notation in hand, we may now introduce our method, called
DOLPHIn (DictiOnary Learning for PHase retrIeval). We consider an
optimization problem which can be interpreted as a combination of an
optimization-based approach to phase retrieval---minimizing the
residual norm with respect to the set of nonlinear equations induced
by the phaseless measurements, cf.~\eqref{eq:WFobjective}---and a
(patch-based) dictionary learning model similar to that used for image
denoising in~\cite{EA06}. The model contains three variables: The
image, or phase retrieval solution $\X$, the dictionary $\D$ and the
matrix $\A$ containing as columns the coefficient vectors of the
representation $\X\approx\cR(\D\A)$. The phase retrieval task consists
of estimating $\X$ and the dictionary learning or sparse coding task
consists of estimating $\D$ and $\A$; a common objective function
provides feedback between the two objectives, with the goal of
improving the phase retrieval reconstruction procedure by encouraging
the patches of~$\X$ to admit a sparse approximation.

Formally, the DOLPHIn formulation consists of minimizing
\begin{align}
  \nonumber\min_{\X,\D,\A}\,&\tfrac14\Norm{\Y-\abs{\cF(\X)}^2}_{\text{F}}^2+\tfrac{\mu}{2}\bignorm{\cE(\X)-\D\A}_{\text{F}}^2+\lambda\sum_{i=1}^{p}\bignorm{\a^i}_1\\
  \label{prob:PRDL}\suchthat~~\,&\X\in[0,1]^{N_1\times N_2},\quad\D\in\cD.
\end{align}
Here, $\norm{\X}_{\text{F}}$ denotes the Frobenius matrix-norm, which
generalizes the Euclidean norm to matrices. The parameters
$\mu,\lambda>0$ in the objective~\eqref{prob:PRDL} provide a way to
control the trade-off between the data fidelity term from the phase
retrieval problem and the approximation sparsity of the image
patches\footnote{We discuss suitable choices and sensitivity of the model to 
these parameters in detail in Section~\ref{subsec:parameters}.}.
To that effect, we use the~$\ell_1$-norm, which is well-known
to have a sparsity-inducing effect~\cite{CDS98}. In order to avoid
scaling ambiguities, we also restrict~$\D$ to be in the subset
$\cD\define\{\D\in\R^{s\times n}\,:\,\norm{\d^j}_2\leq 1~\forall j=1,\dots,n\}$ of matrices with 
column $\l_2$-norms at most $1$, 
and assume $n<p$           
(otherwise, each patch is trivially representable by a $1$-sparse
vector $\a^i$ by including $\x^i/\norm{\x^i}_{2}$ as a column of
$\D$).

The model~\eqref{prob:PRDL} could also easily be modified to include
further side constraints, a different type of nonlinear measurements,
or multiple images or measurements, respectively; we omit these
extensions for simplicity.

\section{Algorithmic Framework}\label{sec:algo}
Similar to classical dictionary learning~\cite{OF96,AEB06,MBP14,T15}
and phase retrieval, problem~(\ref{prob:PRDL}) is nonconvex and
difficult to solve. Therefore, we adopt an algorithm that
provides monotonic decrease of the objective while converging to a
stationary point (see Section~\ref{subsec:monotone} below).

The algorithmic framework we employ is that of \emph{alternating
  minimization}: For each variable $\A$, $\X$ and $\D$ in turn, we
take one step towards solving~\eqref{prob:PRDL} with respect to this
variable alone, keeping the other ones fixed. Each of these
subproblems is convex in the remaining unfixed optimization variable,
and well-known efficient algorithms can be employed accordingly. We
summarize our method in Algorithm~\ref{alg:PRDL}, where the
superscript~$^*$ denotes the adjoint operator (for a matrix $\mZ$,
$\mZ^*$ is thus the conjugate transpose), $\Re(\cdot)$ extracts the
real part of a complex-valued argument, and $\odot$ denotes the
Hadamard (element-wise) product of two matrices. The algorithm also
involves the classical soft-thresholding operator
$\cS_{\tau}(\mZ)\define\max\{0,\abs{\mZ}-\tau\}\odot\sign(\mZ)$ and
the Euclidean projection
$\proj_{\cX}(\mZ)\define\max\{0,\min\{1,\mZ\}\}$ onto $\cX\define
[0,1]^{N_1\times N_2}$; here, all operations are meant component-wise.

To avoid training the dictionary on potentially useless early
estimates, the algorithm performs two phases---while the iteration
counter~$\ell$ is smaller than~$K_1$, the dictionary is not updated.
Below, we explain the algorithmic steps in more detail.

\begin{algorithm*}[tb]
\begin{algorithmic}[1]
  \REQUIRE{Initial image estimate $\X_{(0)}\in[0,1]^{N_1\times N_2}$, initial dictionary $\D_{(0)}\in\cD\subset\R^{s\times n}$, parameters $\mu,\lambda >0$, maximum number of iterations $K_1,K_2$}
  \ENSURE{Learned dictionary $\D=\D_{(K)}$, patch representations $\A=\A_{(K)}$, image reconstructions $\X=\X_{(K)}$ and $\cR(\D\A)$}
  \FOR{$\l=0,1,2,\dots,K_1+K_2\eqqcolon K$}
  \STATE \hspace{-2pt}choose step size $\gamma_{\l}^A$ as explained in Section~\ref{subsec:Aupdate} and update
     $$
     \A_{(\l+1)}\leftarrow \cS_{\gamma_{\l}^{A}\lambda/\mu}\Big(\A_{\l}-\gamma_{\l}^{A}\D_{(\l)}^\top\big(\D_{(\l)}\A_{(\l)}-\cE(\X_{(\l)})\big)\Big)
     $$\label{step:PRDL_Aupdate}
     \STATE \hspace{-2pt}choose step size $\gamma_{\l}^X$ as explained in Section~\ref{subsec:monotone} or \ref{subsec:impl} and update
     $$
     \X_{(\l+1)}\leftarrow \proj_{\cX}\Big(\X_{(\l)}-\gamma_{\l}^X\Big(\Re\Big(\cF^*\big(\cF(\X)\odot(\abs{\cF(\X)}^2-\Y)\big)\Big)
     +\mu\mR\odot\cR\big(\cE(\X)-\D\A\big)\Big)\Big),
     $$\label{step:PRDL_Xupdate}
     \hspace{-0.5em}where~$\mR$ is defined in Section~\ref{subsec:Xupdate}
     \IF{$\l< K_1$}\label{step:PRDL_DupdateStart}
        \STATE do not update the dictionary: $\D_{(\l+1)}\leftarrow\D_{(\l)}$
     \ELSE
         \STATE set $\B\leftarrow\cE(\X_{(\l)})\A_{(\l)}^\top$ and $\mC\leftarrow\A_{(\l)}\A_{(\l)}^\top$
         \FOR{$j=1,\dots,n$}\label{step:PRDL_DupdateLoopstart}
            \IF{$C_{jj}>0$}
               \STATE update $j$-th column: $(\D_{(\l+1)})_{\cdot j}\leftarrow\tfrac{1}{C_{jj}}\big(\B_{\cdot j}-\D_{(\l)}\mC_{\cdot j}\big)+(\D_{(\l)})_{\cdot j}$
            \ELSE
               \STATE reset $j$-th column: e.g., $(\D_{(\l+1)})_{\cdot j}\leftarrow$ random $\cN(0,1)$ vector (in $\R^s$)\label{step:PRDL_DupdateReset}
            \ENDIF
            \STATE project $(\D_{(\l+1)})_{\cdot j}\leftarrow\tfrac{1}{\max\{1,\norm{(\D_{(\l+1)})_{\cdot j}}_2\}}(\D_{(\l+1)})_{\cdot j}$\label{step:PRDL_DupdateEnd}
         \ENDFOR
     \ENDIF
  \ENDFOR
\end{algorithmic}
\caption{Dictionary learning for phase retrieval (DOLPHIn)}
\label{alg:PRDL}
\end{algorithm*}

Note that DOLPHIn actually produces two distinct reconstructions of
the desired image, namely $\X$ (the per se ``image variable'') and
$\cR(\D\A)$ (the image assembled from the sparsely coded
patches)\footnote{Technically, $\cR(\D\A)$ might contain entries not
  in $\cX$, so one should project once more. Throughout, we often omit
  this step for simplicity; differences (if any) between $\cR(\D\A)$
  and $\proj_{\cX}(\cR(\D\A))$ were insignificant in all our tests.}.
Our numerical experiments in Section~\ref{sec:experiments} show that
in many cases, $\cR(\D\A)$ is in fact slightly or even significantly
better than $\X$ with respect to at least one quantitative quality measure and
is therefore also considered a possible reconstruction output of
Algorithm~\ref{alg:PRDL} (at least in the noisy setups we consider in
this paper). Nevertheless, $\X$ is sometimes more visually appealing
 and can be used, for instance, to refine parameter
settings (if it differs strongly from $\cR(\D\A)$) or to assess the
influence of the patch-based ``regularization'' on the pure non-sparse
Wirtinger Flow method corresponding to the formulation where~$\lambda$
and~$\mu$ are set to zero.

\subsection{Updating the Patch Representation Vectors}\label{subsec:Aupdate}
Updating $\A$ (i.e., considering \eqref{prob:PRDL} with $\D$ and $\X$ fixed at
their current values) consists of decreasing the objective
\begin{equation}\label{prob:Aupdate}
  \sum_{i=1}^{p}\left(\tfrac{1}{2}\bignorm{\D_{(\l)}\a^i -\x^{i}_{(\l)}}_2^2 +\tfrac{\lambda}{\mu}\bignorm{\a^i}_1\right),
\end{equation}
which is separable in the patches~$i=1\,\ldots,p$. Therefore, we can update all
vectors~$\a^i$ independently and/or in parallel. To do so, we choose to
perform one step of the well-known algorithm ISTA (see,~e.g.,
\cite{BT09}), which
is a gradient-based method that is able to take into account a non-smooth
regularizer such as the~$\ell_1$-norm.
Concretely, the following update is performed for each
$i=1,\dots,p$:
\begin{equation}\label{eq:colAupdate}
  \a^i_{(\l+1)} = \cS_{\gamma_{\l}^{A}\lambda/\mu}\Big(\a^i_{(\l)}-\gamma_{\l}^{A}
  \D^\top_{(\l)}\big(\D_{(\l)}\a^i_{(\l)}-\x^i_{(\l)}\big)\Big).
\end{equation}
This update involves a gradient descent step (the gradient with
respect to $\a^i$ of the smooth term in each summand
of~\eqref{prob:Aupdate} is
$\D^\top_{(\l)}\big(\D_{(\l)}\a^i_{(\l)}-\x^i_{\l}\big)$,
respectively) followed by soft-thresholding. Constructing
$\A_{(\l+1)}$ from the $\a^i_{(\l+1)}$ as specified above is
equivalent to Step~\ref{step:PRDL_Aupdate} of
Algorithm~\ref{alg:PRDL}.

The step size parameter $\gamma_{\l}^A$ can be chosen in $(0,1/L_A)$,
where $L_A$ is an upper bound on the Lipschitz constant of the
gradient; here, $L_A=
\norm{\D_{(\l)}^\top\D_{(\l)}}_{2}=\norm{\D_{(\l)}}_2^2$ would be
appropriate, but a less computationally demanding strategy is to use a
backtracking scheme to automatically update~$L_A$~\cite{BT09}. 

A technical subtlety is noteworthy in this regard: We can either find
one $\gamma_{\l}^A$ that works for the whole matrix-variable update
problem---this is what is stated implicitly in
Step~\ref{step:PRDL_Aupdate}---or we could find different values, say
$\gamma_{\l}^{a,i}$, for each column $\a^i$, $i=1,\dots,p$, of $\A$
separately. Our implementation does the latter, since it
employs a backtracking strategy for each column update independently.

\subsection{Updating the  Image Estimate}\label{subsec:Xupdate}
With $\D=\D_{(\l)}$ and $\A=\A_{(\l+1)}$ fixed, updating
$\X$ consists of decreasing the objective
\begin{align}\label{prob:Xupdate}
  &\tfrac14\bignorm{\Y-\abs{\cF(\X)}^2}_F^2+\tfrac{\mu}{2}\bignorm{\cE(\X)-\D\A}_F^2\\
   \nonumber\text{with}\quad&\X\in\cX=[0,1]^{N_1\times N_2}.
\end{align}
This problem can be seen as a regularized version of the phase
retrieval problem (with regularization parameter~$\mu$) that
encourages the patches of~$\X$ to be close to the sparse
approximation~$\D\A$ obtained during the previous (inner) iterations.

Our approach to decrease the value of the objective
\eqref{prob:Xupdate} is by a projected gradient descent step. In fact,
for $\mu=0$, this step reduces to the Wirtinger flow method
\cite{CLS14}, but with necessary modifications to take into account
the constraints on $\X$ (real-valuedness and variable bounds $[0,1]$).

The gradient of
$\phi(\X)\define\tfrac14\norm{\Y-\abs{\cF(\X)}^2}_F^2$ with respect to
$\X$ can be computed as
\[
  \nabla\phi(\X)=\Re\Big(\cF^*\big(\cF(\X)\odot(\abs{\cF(\X)}^2-\Y)\big)\Big),
\]
by using the chain rule. For
$\psi(\X)\define\tfrac{\mu}{2}\norm{\cE(\X)-\D\A}_F^2$, the gradient is given
by
\[
   \nabla\psi(\X)=\mu\cE^*\big(\cE(\X)-\D\A\big)=\mu\mR\odot\cR\big(\cE(\X)-\D\A\big),
\]
where $\mR$ is an $N_1 \times N_2$ matrix whose entries $r_{ij}$ equal
the number of patches the respective pixel $x_{ij}$ is contained
in. Note that if the whole image is divided into a complete set of
nonoverlapping patches, $\mR$ will just be the all-ones matrix;
otherwise, the element-wise multiplication with $\mR$ undoes the
averaging of pixel values performed by $\cR$ when assembling an image
from overlapping patches.

Finally, the gradient w.r.t. $\X$ of the objective
in~\eqref{prob:Xupdate} is
$\nabla\phi(\X)+\nabla\psi(\X)\in\R^{N_1\times N_2}$, and the update
in Step~\ref{step:PRDL_Xupdate} of Algorithm~\ref{alg:PRDL} is indeed
shown to be a projected gradient descent step. Typically, a
backtracking (line search) strategy is used for choosing the step
size~$\gamma_{\ell}^X$; see Theorem~\ref{thm:convCGD} in Section~\ref{subsec:monotone} for a selection
 rule that gives theoretical convergence, and also 
Section~\ref{subsec:impl} for a heuristic alternative.

\subsection{Updating the Dictionary}\label{subsec:Dupdate}
To update the dictionary, i.e., to approximately
solve~\eqref{prob:PRDL} w.r.t. $\D$ alone, keeping $\X$ and $\A$ fixed
at their current values, we employ one pass of a block-coordinate descent (BCD) algorithm
on the columns of the dictionary~\cite{MBPS10}.
The objective to decrease may be written as
\begin{equation}\label{prob:Dupdate}
  \tfrac12\sum_{i=1}^{p}\bignorm{\D\a^i_{(\l+1)} -\x^i_{(\l+1)}}_2^2\quad\suchthat\quad\D\in\cD,
\end{equation}
and the update rule given by
Steps~\ref{step:PRDL_DupdateStart}\,--\ref{step:PRDL_DupdateEnd}
corresponds\footnote{In \cite[Algo.~11]{MBP14}, and in our
  implementation, we simply normalize the columns of $\D$; it is
  easily seen that any solution with $\norm{\d^j}_2<1$ for some $j$ is
  suboptimal (w.r.t.~\eqref{prob:PRDL}) since raising it to $1$ allows
  to reduce coefficients in $\A$ and thus to improve the $\l_1$-term
  of the DOLPHIn objective~\eqref{prob:PRDL}. However, using the
  projection is more convenient for proving the convergence results
  without adding more technical subtleties w.r.t. this aspect.} to one
iteration of \cite[Algorithm~11]{MBP14} applied
to~\eqref{prob:Dupdate}. 

To see this, note that each column update problem has a
closed-form solution:
\begin{align*}
  (\d^j)_{(\l+1)}= \proj_{\norm{\cdot}_2\leq 1}\left(\frac{1}{\sum_{i=1}^{p}(a^i_j)^2}\sum_{i=1}^p a^i_j\big(\x^i-\sum_{\begin{subarray}{c}k=1\\k\neq j\end{subarray}}^{n}a^i_k\d^k\big)\right)
  =\frac{1}{\max\{1,\norm{\tfrac{1}{w_j}\q^j}_2\}}\big(\tfrac{1}{w_j}\q^j\big)
\end{align*}
with $w_j\define\sum_{i} (a^i_j)^2$ and $\q^j\define\sum_{i}
a^i_j\big(\x^i-\sum_{k\neq j}a^i_k\d^k\big)$; here, we abbreviated
$\a^i\define\a^i_{(\ell+1)}$, $\x^i\define\x^i_{(\ell+1)}$. If
$w_j=0$, and thus $a^i_j=0$ for all $i=1,\dots,p$, then column $\d^j$
is not used in any of the current patch representations; in that case,
the column's update problem has a constant objective and is therefore
solved by any $\d$ with $\norm{\d}_2\leq 1$, e.g., a
normalized random vector as in Step~12 of Algorithm~1.
The computations performed in
Steps~\ref{step:PRDL_DupdateLoopstart}--\ref{step:PRDL_DupdateEnd} of
Algorithm~\ref{alg:PRDL} are equivalent to these solutions,
expressed differently using the matrices $\B$ and
$\mC$ defined there. Note that the operations
 could be parallelized to speed up computation.

\subsection{Convergence of the Algorithm}\label{subsec:monotone}
As mentioned at the beginning of this section, with appropriate step
size choices, DOLPHIn (Algorithm~\ref{alg:PRDL}) exhibits the property
of monotonically decreasing the objective function
value~\eqref{prob:PRDL} at each iteration. In particular, many
line-search type step size selection mechanisms aim precisely at
reducing the objective; for simplicity, we will simply refer to such
subroutines as ``suitable backtracking schemes'' below. Concrete
examples are the ISTA backtracking from~\cite[Section~3]{BT09} we can
employ in the update of $\A$, or the rule given in
Theorem~\ref{thm:convCGD} for the projected gradient descent update of
$\X$ (a different choice is described in Section~\ref{subsec:impl});
further variants are discussed, e.g., in~\cite{I03}.

\begin{prop}\label{prop:monotone}
  Let $(\A_{(\l)},\X_{(\l)},\D_{(\l)})$ be the current iterates (after
  the $\l$-th iteration) of Algorithm~\ref{alg:PRDL} with step
  sizes $\gamma_{\ell}^{X}$ and $\gamma_{\ell}^{A}$ determined by
  suitable backtracking schemes (or arbitrary $0<\gamma_{\l}^{A}<
  1/\norm{\D_{(\l)}^\top\D_{(\l)}}_2$, resp.) and let $f_{i,j,k}$
  denote the objective function value of the DOLPHIn
  model~\eqref{prob:PRDL} at $(\A_{(i)},\X_{(j)},\D_{(k)})$. Then,
  DOLPHIn either terminates in the $(\l+1)$-th iteration, or
  it holds that $f_{\l+1,\l+1,\l+1}\leq f_{\l,\l,\l}$.
\end{prop}
\begin{IEEEproof}
  Since we use ISTA to update $\A$, it follows from~\cite{BT09} that 
	$f_{\l+1,\l,\l}\leq f_{\l,\l,\l}$. Similarly, a
  suitable backtracking strategy is known to enforce descent in the
  projected gradient method when the gradient is locally Lipschitz-continuous, whence $f_{\l+1,\l+1,\l}\leq f_{\l+1,\l,\l}$. Finally,
  $f_{\l+1,\l+1,\l+1}\leq f_{\l+1,\l+1,\l}$ follows from standard
  results for BCD methods applied to convex problems, see, e.g.,
  \cite{B99}. Combining these inequalities proves the claim.
\end{IEEEproof}

  The case of termination in Proposition~\ref{prop:monotone} can occur when
  the backtracking scheme is combined with a maximal number of trial
  steps, which are often used as a safeguard against numerical
  stability problems or as a heuristic stopping condition to terminate
  the algorithm if no (sufficient) improvement can be reached even
  with tiny step sizes.  Note also that the assertion of
  Proposition~\ref{prop:monotone} trivially holds true if all step sizes are
  $0$; naturally, a true descent of the objective requires a strictly
  positive step size in at least one update step. 
  In our algorithm, step size positivity can always be guaranteed since all
  these updates involve objective functions whose gradient is Lipschitz
  continuous, and backtracking essentially finds step sizes inversely
  proportional to (local) Lipschitz constants. (Due to non-expansiveness, the
  projection in the $\X$-update poses no problem either).
\smallskip

Adopting a specific Armijo-type step size selection rule
for the $\X$-update allows us to infer a convergence result, stated 
in Theorem~\ref{thm:convCGD} below. To
simplify the presentation, let
$f^X_{\l}(\X)\define\tfrac14\norm{\Y-\abs{\cF(\X)}^2}_{\text{F}}^2+\tfrac{\mu}{2}\norm{\cE(\X)-\D_{(\l)}\A_{(\l+1)}}_{\text{F}}^2$
and ${\bf\Gamma}^X_{\l}\define\nabla f^X_{\l}(\X_{(\l)})$
(cf. $\nabla\phi(\X)+\nabla\psi(\X)$ in Section~\ref{subsec:Xupdate}).

\begin{thm}\label{thm:convCGD}
  Let $\eta\in(0,1)$ and $\bar\eta>0$. Consider the DOLPHIn variant
  consisting of Algorithm~\ref{alg:PRDL} with $K_2=\infty$ and the
  following Armijo rule to be used in Step~\ref{step:PRDL_Xupdate}:\\[-0.5em]
  \begin{center}
    \begin{minipage}{0.9\textwidth}
      Determine $\gamma^X_{\l}$ as the largest number
      in $\{\bar\eta\eta^k\}_{k=0,1,2,\dots}$ such that
      $f^X_{\l}(\proj_{\cX}(\X_{(\l)}-\gamma^X_{\l} {\bf\Gamma}^X_{\l}))-f^X_{\l}(\X_{(\l)})\leq 
      -\tfrac{1}{2\gamma^X_{\l}} \norm{\proj_{\cX}(\X_{(\l)}-\gamma^X_{\l} {\bf\Gamma}^X_{\l})-\X_{(\l)}}_{\text{F}}^2$
      and set $\X_{(\l+1)}\define\proj_{\cX}(\X_{(\l)}-\gamma^X_{\l} {\bf\Gamma}^X_{\l})$.
    \end{minipage}
  \end{center}
  ~\\\noindent If $\mu\geq 1$ and,  for some $0<\underline{\nu}\leq\bar\nu$, it holds that
  $\underline{\nu}\leq\gamma^X_{\l},\gamma^A_{\l}$ (or
  $\gamma^{a,i}_{\l}$), $\sum_{i=1}^{p}(\a^i_{(\l)})_j^2\leq\bar\nu$
  for all $\l$ and $j$ (and $i$), then
  every accumulation point of the sequence
  $\{(\A_{(\l)},\X_{(\l)},\D_{(\l)})\}_{\l=0,1,2,\dots}$ of DOLPHIn
  iterates is a stationary point of problem~\eqref{prob:PRDL}.
\end{thm}
\begin{IEEEproof}
 The proof works by expressing Algorithm~\ref{alg:PRDL} as a specific
  instantiation of the coordinate gradient descent (CGD) method
  from~\cite{TY09} and analyzing the objective descent achievable in
  each update of $\A$, $\X$ and $\D$, respectively. The technical details make the rigorous formal proof somewhat lengthy; therefore, we only sketch it here
  and defer the full proof to Appendix~\ref{sec:appendixB}.

  The CGD method works by solving subproblems to obtain directions of
  improvement for blocks of variables at a time---in our case, (the
  columns of) $\A$, the matrix $\X$, and the columns of $\D$
  correspond to such blocks---and then taking steps along these
  directions. More specifically, the directions are generated using a
  (suitably parameterized) strictly convex quadratic approximation of
  the objective (built using the gradient). Essentially due to the
  strict convexity, it is then always possible to make a
  \emph{positive} step along such a direction that decreases the
  (original) objective, unless stationarity already holds. Using a
  certain Armijo line-search rule designed to find such positive step
  sizes which achieve a sufficient objective reduction,
  \cite[Theorem~1]{TY09} ensures (under mild further assumptions,
  which in our case essentially translate to the stated boundedness
  requirement of the step sizes) that every accumulation point of the
  iterate sequence is indeed a stationary point of the addressed
  (block-separable) problem.

  To embed DOLPHIn into the CGD framework, we can interpret the
  difference between one iterate and the next (w.r.t. the variable
  ``block'' under consideration) as the improvement direction, and
  proceed to show that we can always choose a step size equal to 1 in
  the Armijo-criterion from~\cite{TY09} (cf. (9) and (46)
  therein). For this to work out, we need to impose slightly stricter
  conditions on other parameters used to define that rule than what is
  needed in \cite{TY09}; these conditions are derived directly from
  known descent properties of the $\D$- and $\A$-updates of our method
  (essentially, ISTA descent properties as in~\cite{BT09}). That way,
  the $\D$- and $\A$-updates automatically satisfy the specific CGD
  Armijo rule, and the actual backtracking scheme for the $\X$-update
  given in the present theorem can be shown to assert that our
  $\X$-update does so as well. (The step sizes used in DOLPHIn could
  also be reinterpreted in the CGD framework as scaling factors of
  diagonal Hessian approximations of the combined objective to be used
  in the direction-finding subproblems. With such simple Hessian
  approximations, the obtained directions are then indeed equivalent
  to the iterate-differences resulting from the DOLPHIn update
  schemes.) The claim then follows directly from
  \cite[Theorem~1(e) (and its extensions discussed in Section~8)]{TY09}.
\end{IEEEproof}

  A more formal explanation for why the step sizes can be
  chosen positive in each step can be found on page~392
  of~\cite{TY09}; the boundedness of approximate Hessians is stated
  in~\cite[Assumption~1]{TY09}. Arguably, assuming step sizes are
  bounded away from zero by a constant may become problematic in
  theory (imagine an Armijo-generated step size sequence converging to
  zero), but will not pose a problem in practice where one always
  faces the limitations of numerical machine precision. (Note
  also that, in practice, the number of line-search trials can be
  effectively reduced by choosing $\overline{\eta}$ based on the
  previous step size~\cite{TY09}.)
	
	Our implementation uses a different backtracking scheme for
  the $\X$-update (see Section~\ref{subsec:impl}) that can be viewed
  as a cheaper heuristic alternative to the stated Armijo-rule which
  still ensures monotonic objective descent (and hence is ``suitable''
  in the context of Proposition~\ref{prop:monotone}), also enables strictly
  positive steps, and empirically performs equally well. Finally, we
  remark that the condition $\mu\geq 1$ in Theorem~\ref{thm:convCGD}
  can be dropped if the relevant objective parts of
  problem~\eqref{prob:PRDL} are not rescaled for the $\A$- and
  $\D$-updates, respectively.
\smallskip

To conclude the discussion of convergence, we point out
that one can obtain a linear rate of convergence for DOLPHIn with the Armijo
rule from Theorem~\ref{thm:convCGD}, by extending the results of
\cite[Theorem~2 (cf. Section~8)]{TY09}.


\section{Numerical Results}\label{sec:experiments}
In this section, we discuss various numerical experiments to study the
effectiveness of the DOLPHIn algorithm. To that end, we consider
several types of linear operators $\cF$ within our
model~\eqref{prob:PRDL} (namely, different types of Gaussian random
operators and coded diffraction models).  Details on the
expe\-ri\-mental setup and our implementation are given in the first
two subsections, before presenting the main numerical results in
Subsection~\ref{subsec:experiments}. Our experiments demonstrate that
with noisy measurements, DOLPHIn gives significantly better image
reconstructions than the Wirtinger Flow method~\cite{CLS14}, one
recent state-of-the-art phase retrieval algorithm, thereby showing
that introducing sparsity via a (simultaneously) learned dictionary is
indeed a promising new approach for signal reconstruction from noisy,
phaseless, nonlinear measurements. Furthermore, we discuss sensitivity
of DOLPHIn with regard to various (hyper-)parameter choices
(Subsections.~\ref{subsec:parameters},~\ref{subsec:inneriters} and~\ref{subsec:firstphase}) and a
variant in which the $\l_1$-regularization term in the objective is
replaced by explicit constraints on the sparsity of the
patch-representation coefficient vectors $\a^i$
(Subsection.~\ref{subsec:l0prdl}).

\subsection{Experimental Setup}\label{subsec:imp}
We consider several linear operators $\cF$ corresponding to different
types of measurements that are classical in the phase retrieval
literature.  We denote by $\F$ the (normalized) 2D-Fourier operator
(implemented using fast Fourier transforms), and introduce two complex
Gaussian matrices $\G\in\C^{M_1\times N_1}$, $\H\in\C^{M_2\times
  N_2}$, whose entries are i.i.d. samples from the distribution
$\cN(0,\I/2)+i\cN(0,\I/2)$.  Then, we experiment with the operators
$\cF(\X)=\G\X$, $\cF(\X)=\G\X\G^*$, $\cF(\X)=\G\X\H^*$, and the coded
diffraction pattern model
\begin{equation}\label{eq:CDPop}
   \cF(\X)=\left(\begin{array}{c}\F\big(\overline{\M}_1\odot\X\big)\\[-0.5em]\vdots\\\F\big(\overline{\M}_m\odot\X\big)\end{array}\right),~~\cF^*(\mZ)=\sum_{j=1}^{m}\big(\M_j\odot\F^*(\mZ_j)\big),
\end{equation}
where $\mZ_j\define\mZ_{\{(j-1)N_{1}+1,\dots,j N_{1}\},\cdot}$ (i.e., $\mZ^\top=(\mZ_1^\top,\dots,\mZ_m^\top)$) and the $\M_j$'s are admissible coded diffraction patterns (CDPs),
see for instance \cite[Section~4.1]{CLS14}; in our experiments we used
ternary CDPs, such that each $\M_j$ is in $\{0,\pm 1\}^{N_1\times
  N_2}$. (Later, we will also consider octanary CDPs with
$\M_j\in\{\pm\sqrt{2}/2,\pm i\sqrt{2}/2,\pm\sqrt{3},\pm i\sqrt{3}\}\in\C^{N_1\times N_2}$.) 

To reconstruct $\hat\X$, we choose an oversampling setting where
$M_1=4N_1$, $M_2=4N_2$ and/or $m=2$, respectively. Moreover, we
corrupt our measurements with additive white Gaussian noise~$\mN$ such
that $\text{SNR}(\Y,\abs{\cF(\hat\X)}^2+\mN)=10$\,dB for the
Gaussian-type, and $20$\,dB for CDP measurements, respectively. Note
that these settings yield, in particular, a relatively heavy noise
level for the Gaussian cases and a relatively low oversampling ratio
for the CDPs.

\subsection{Implementation Details}\label{subsec:impl}
We choose to initialize our algorithm with a simple random image
$\X_{(0)}$ in $\cX$ to demonstrate the robustness of our approach with
respect to its initialization. Nevertheless, other choices are
possible. For instance, one may also initialize $\X_{(0)}$ with a
power-method scheme similar to that proposed in~\cite{CLS14}, modified
to account for the real-valuedness and box-constraints. The dictionary
is initialized as $\D_{(0)}=(\I,\F_D)$ in $\R^{s\times 2s}$, where
$\F_D$ corresponds to the two-dimensional discrete cosine transform
(see, e.g.,~\cite{EA06}).

To update $\A$, we use the ISTA implementation from the SPAMS
package\footnote{\url{http://spams-devel.gforge.inria.fr/}}~\cite{MBPS10}
with its integrated backtracking line search (for $L_A$).  Regarding
the step sizes $\gamma_{\l}^X$ for the update of~$\X$
(Step~\ref{step:PRDL_Xupdate} of Algorithm~\ref{alg:PRDL}), we adopt the
following simple strategy, which is similar to that from~\cite{BT09}
and may be viewed as a heuristic to the Armijo rule from
Theorem~\ref{thm:convCGD}: Whenever
the gradient step leads to a reduction in the objective function
value, we accept it. Otherwise, we recompute the step with
$\gamma_{\l}^X$ halved until a reduction is achieved; here, as a
safeguard against numerical issues, we implemented a limit of 100
trials (forcing termination in case all were unsuccessful), but this
was never reached in any of our computational experiments. Regardless
of whether $\gamma_{\l}^X$ was reduced or not, we reset its value to
$1.68\gamma_{\l}^X$ for the next round; the initial step
size 
is~$\gamma_0^X=10^4/f_{(0)}$, where $f_{(0)}$ is the objective
function of the DOLPHIn model~\eqref{prob:PRDL}, evaluated at
$\X_{(0)}$, $\D_{(0)}$ and least-squares patch representations
$\arg\min_{\A}\norm{\cE(\X_{(0)})-\D_{(0)}\A}_{\text{F}}^2$. (Note
that, while this rule deviates from the theoretical convergence
Theorem~\ref{thm:convCGD}, Propositon~\ref{prop:monotone} and
the remarks following it remain applicable.)

Finally, we consider nonoverlapping $8\times 8$ patches and run
DOLPHIn (Algorithm~\ref{alg:PRDL}) with $K_1=25$ and $K_2=50$; the
regularization/penalty parameter values can be read from
Table~\ref{tab:PRDLresults} (there, $m_Y$ is the number of elements
of~$\Y$). We remark that these parameter values were empirically
benchmarked to work well for the measurement setups and instances
considered here; a discussion about the stability of our approach with
respect to these parameter choices is presented below in
Section~\ref{subsec:parameters}. Further experiments with a
sparsity-constrained DOLPHIn variant and using overlapping patches are
discussed in Section~\ref{subsec:l0prdl}.

Our DOLPHIn code is available online on the first author's webpage\footnote{\url{http://www.mathematik.tu-darmstadt.de/~tillmann/}}.

\subsection{Computational Experiments}\label{subsec:experiments}
\begin{figure*}[tb]
\centering
\subfigure[]{\includegraphics[width=0.23\textwidth,resolution=300]{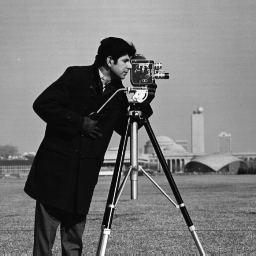}}
~
\subfigure[]{\includegraphics[width=0.23\textwidth,resolution=300]{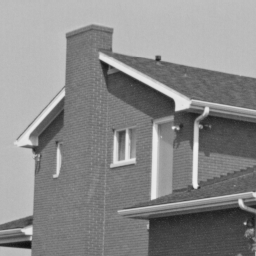}}
~
\subfigure[]{\includegraphics[width=0.23\textwidth,resolution=300]{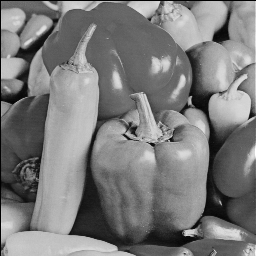}}
~
\subfigure[]{\includegraphics[width=0.23\textwidth,resolution=300]{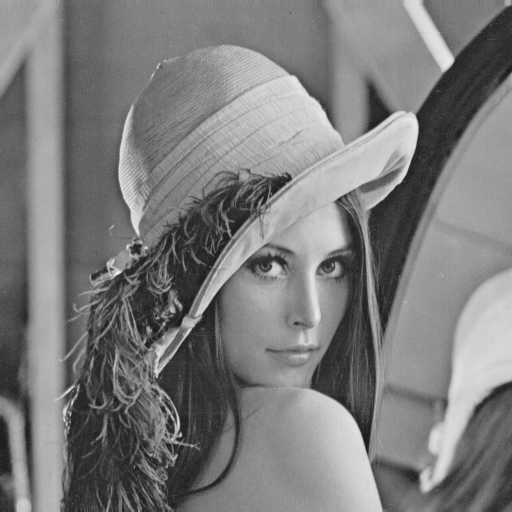}}
~
\subfigure[]{\includegraphics[width=0.23\textwidth,resolution=300]{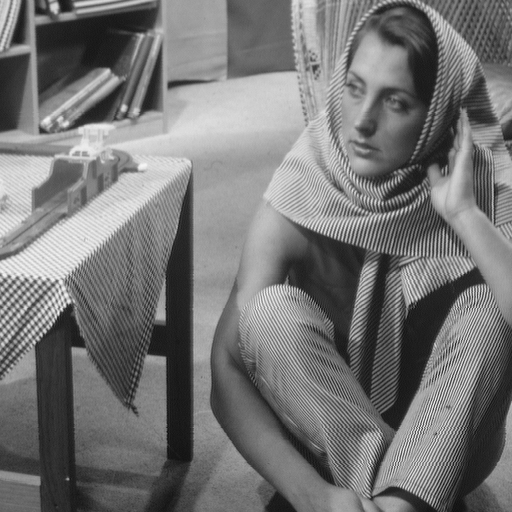}}
~
\subfigure[]{\includegraphics[width=0.23\textwidth,resolution=300]{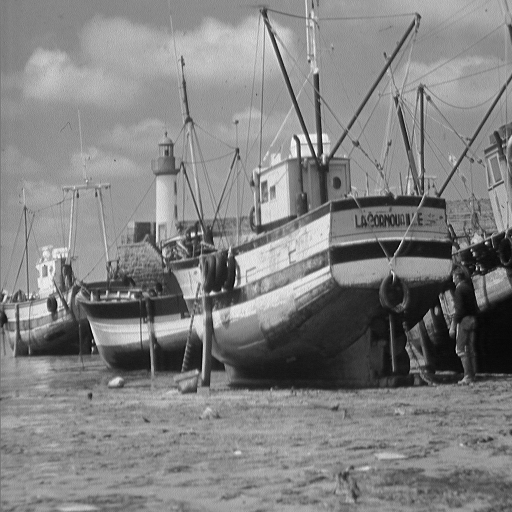}}
~
\subfigure[]{\includegraphics[width=0.23\textwidth,resolution=300]{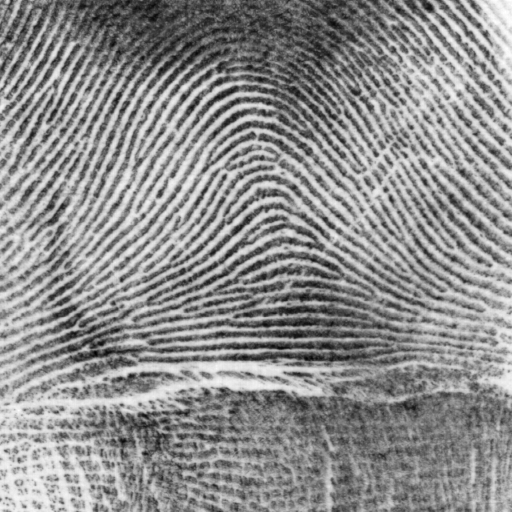}}
~
\subfigure[]{\includegraphics[width=0.23\textwidth,resolution=300]{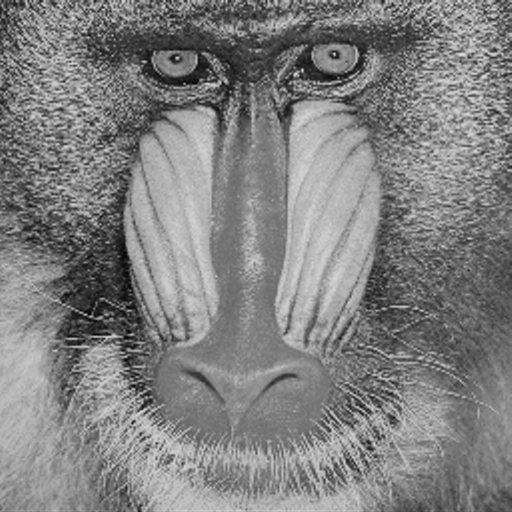}}

\caption{Test images. (a)--(c): cameraman, house and peppers (size
$256\times 256$). (d)--(h): lena, barbara, boat, fingerprint and mandrill (size $512\times 512$).}
\label{fig:images}
\end{figure*}

We test our method on a collection of typical (grayscale) test images
used in the literature, see Figure~\ref{fig:images}. All experiments
were carried out on a Linux $64$-bit quad-core machine ($2.8$\,GHz,
$8$\,GB RAM) running Matlab R2016a (single-thread).

\begin{table*}[tb]
  \caption{Test results for $m_Y$ Gaussian-type and coded diffraction pattern (CDP) measurements. We report mean values (geometric mean for CPU times) per measurement type, obtained from three instances with random $\X_{(0)}$ and random noise for each of the three $256\times 256$ and five $512\times 512$ images, w.r.t. the reconstructions from DOLPHIn ($\X_{\text{DOLPHIn}}$ and $\proj_{\cX}(\cR(\D\A))$) with parameters $(\mu,\lambda)$ and (real-valued, $[0,1]$-constrained) Wirtinger Flow ($\X_{\text{WF}}$), respectively. (CPU times in seconds, PSNR in decibels).}
  \label{tab:PRDLresults}
  \centering
	\begin{tabular*}{\textwidth}{@{\extracolsep{\fill}}l@{\quad}l@{\quad}c@{\quad}c@{\quad}c@{\quad}c@{\quad}c@{\quad}c@{\quad}c@{\quad}c@{\quad}c@{\quad}c@{\quad}c@{\quad}}
\toprule
 & & \multicolumn{5}{c}{$256\times 256$ instances} & \multicolumn{5}{c}{$512\times 512$ instances}\\\cmidrule(r{9pt}){3-7}\cmidrule(r{9pt}){8-12}
$\cF$ type & reconstruction & $(\mu,\lambda)/m_Y$ & time & PSNR & SSIM & $\varnothing$\,$\norm{\a^i}_0$ & $(\mu,\lambda)/m_Y$ & time & PSNR & SSIM & $\varnothing$\,$\norm{\a^i}_0$\\
\midrule
$\G\hat\X$ & $\X_{\text{DOLPHIn}}$ & (0.5,0.105) & 12.69 & 24.69 & 0.5747 & & (0.5,0.105) & 68.62 & 24.42 & 0.6547 & \\
& $\proj_{\cX}(\cR(\D\A))$ & & & 23.08 & 0.6644 & 3.77 & & & 22.66 & 0.6807 & 6.30\\
& $\X_{\text{WF}}$ & & 7.49 & 19.00 & 0.2898 & -- & & 49.23 & 18.83 & 0.3777 & --\\[0.5em]
$\G\hat\X\G^*$ & $\X_{\text{DOLPHIn}}$ & (0.5,0.210) & 51.78 & 22.67 & 0.4135 & & (0.5,0.210) & 357.49 & 22.59 & 0.5296 & \\
& $\proj_{\cX}(\cR(\D\A))$ & & & 23.70 & 0.7309 & 7.45 & & & 23.43 & 0.7685 & 11.37\\
& $\X_{\text{WF}}$ & & 47.76 & 22.66 & 0.4131 & -- & & 349.28 & 22.58 & 0.5290 & --\\[0.5em]
$\G\hat\X\H^*$ & $\X_{\text{DOLPHIn}}$ & (0.5,0.210) & 52.18 & 22.67 & 0.4132 & & (0.5,0.210) & 357.66 & 22.57 & 0.5286 & \\
& $\proj_{\cX}(\cR(\D\A))$ & & & 23.68 & 0.7315 & 7.50 & & & 23.43 & 0.7667 & 11.38\\
& $\X_{\text{WF}}$ & & 48.24 & 22.65 & 0.4127 & -- & & 348.54 & 22.55 & 0.5282 & --\\[0.5em]
CDP (cf.~\eqref{eq:CDPop}) & $\X_{\text{DOLPHIn}}$ & (0.05,0.003) & 8.56 & 27.15 & 0.7416 & & (0.05,0.003) & 36.72 & 27.33 & 0.7819 & \\
& $\proj_{\cX}(\cR(\D\A))$ & & & 26.58 & 0.7654 & 7.85 & & & 26.33 & 0.7664 & 11.48\\
& $\X_{\text{WF}}$ & & 2.83 & 13.10 & 0.1170 & -- & & 14.79 & 12.70 & 0.1447 & --\\
\bottomrule
\end{tabular*}
\end{table*}

We evaluate our approach with the following question in mind:
\emph{Can we improve upon the quality of reconstruction compared to
  standard phase retrieval algorithms?}  Standard methods cannot
exploit sparsity if the underlying basis or dictionary is unknown; as
we will see, the introduced (patch-) sparsity indeed allows for better
recovery results (at least in the oversampling and noise regimes
considered here).

To evaluate the achievable sparsity, we look at the average number of
nonzeros in the columns of $\A$ after running our
algorithm. Generally, smaller values indicate an improved suitability
of the learned dictionary for sparse patch coding (high values often
occur if the regularization parameter $\lambda$ is too small and the
dictionary is learning the noise, which is something we would like to
avoid). To assess the quality of the image reconstructions, we
consider two standard measures, namely the peak signal-to-noise ratio
(PSNR) of a reconstruction as well as its structural similarity index
(SSIM) \cite{WBSS04}. For PSNR, larger values are better, and
SSIM-values closer to~$1$ (always ranging between $0$ and $1$)
indicate better visual quality.

\begin{figure*}[tb]
\centering
\subfigure[]{\includegraphics[width=0.23\textwidth,height=0.23\textwidth,resolution=300]{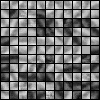}} 
~
\subfigure[]{\includegraphics[width=0.23\textwidth,resolution=300]{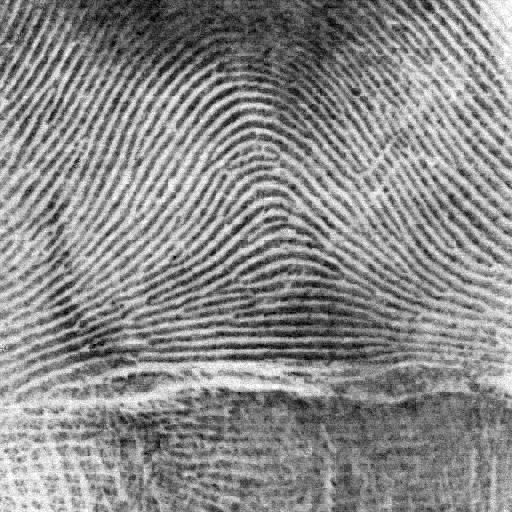}}
~
\subfigure[]{\includegraphics[width=0.23\textwidth,resolution=300]{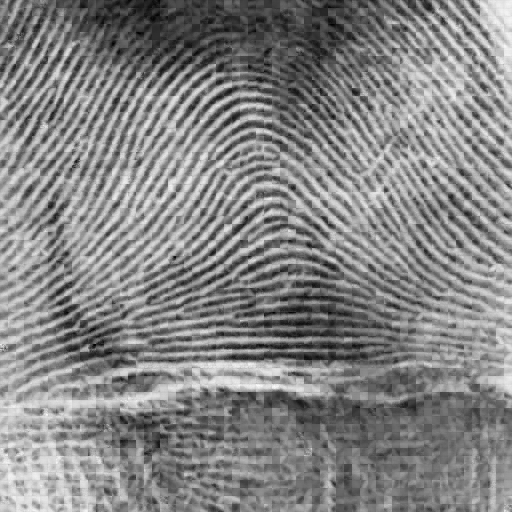}}
~
\subfigure[]{\includegraphics[width=0.23\textwidth,resolution=300]{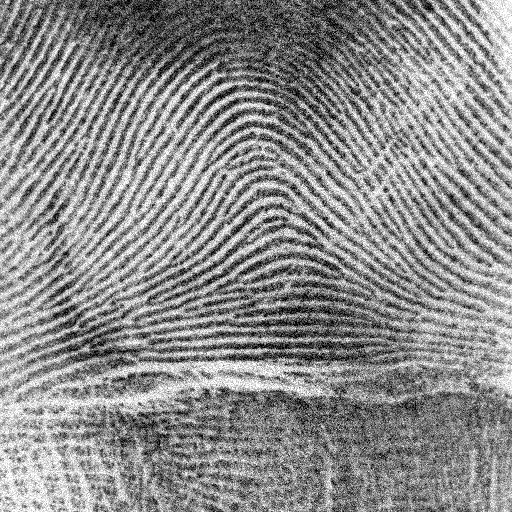}} 
\caption{DOLPHIn example: Image original is the $512\times 512$
  ``fingerprint'' picture, measurements are noisy Gaussian
  $\G\hat{\X}$ ($M_1=4N_1$, noise-SNR $10$\,dB),
  $(\mu,\lambda)=(0.5,0.105)m_Y$. (a) final dictionary (excerpt, magnified), (b)
  image reconstruction $\X_{\text{DOLPHIn}}$, (c) image reconstruction
  $\cR(\D\A)$ from sparsely coded patches, (d) reconstruction
  $\X_{\text{WF}}$ after $75$ WF iterations. Final PSNR values:
  $22.37$\,dB for $\cR(\D\A)$, $23.74$\,dB for $\X_{\text{DOLPHIn}}$,
  $18.19$\,dB for $\X_{\text{WF}}$; final SSIM values: $0.7903$ for
  $\cR(\D\A)$, $0.8152$ for $\X_{\text{DOLPHIn}}$, $0.5924$ for
  $\X_{\text{WF}}$; average $\norm{\a^i}_0$ is $10.06$.}
\label{fig:PRDLexample2}
\end{figure*}

\begin{figure*}[tb]
\centering
\subfigure[]{\includegraphics[width=0.32\textwidth,resolution=300]{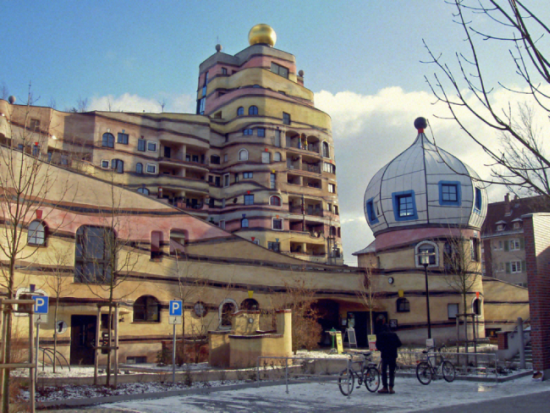}} 
~
\subfigure[]{\includegraphics[width=0.32\textwidth,resolution=300]{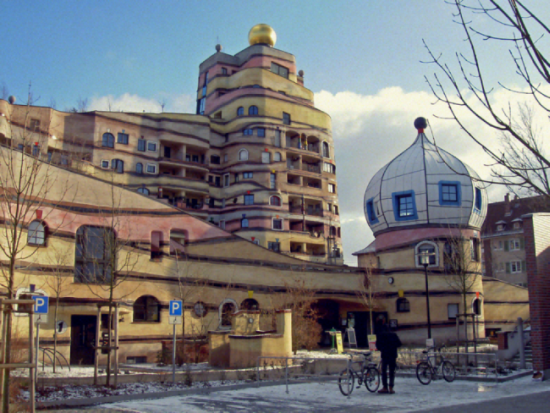}}
~
\subfigure[]{\includegraphics[width=0.32\textwidth,resolution=300]{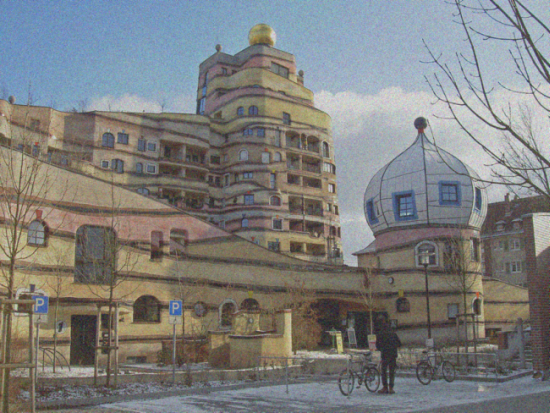}}
\caption{DOLPHIn example: Image original is a $2816\times 2112$ photo
  of the ``Waldspirale'' building in Darmstadt, Germany; measurements
  are noisy CDPs (obtained using two ternary masks), noise-SNR $20$\,dB, 
  $(\mu,\lambda)=(0.05,0.007)m_Y$. (a) image reconstruction
  $\X_{\text{DOLPHIn}}$, (b) image reconstruction $\cR(\D\A)$ from
  sparsely coded patches, (c) reconstruction $\X_{\text{WF}}$ after
  $75$ WF iterations. Final PSNR values: $23.40$\,dB for $\cR(\D\A)$,
  $24.72$\,dB for $\X_{\text{DOLPHIn}}$, $12.63$\,dB for
  $\X_{\text{WF}}$; final SSIM values: $0.6675$ for $\cR(\D\A)$,
  $0.6071$ for $\X_{\text{DOLPHIn}}$, $0.0986$ for $\X_{\text{WF}}$;
  average $\norm{\a^i}_0$ is $12.82$. (Total reconstruction time
  roughly $30$\,min (DOLPHIn) and $20$\,min (WF), resp.) \emph{Original image taken from Wikimedia Commons, under Creative Commons Attribution-Share Alike 3.0 Unported license.}}
\label{fig:PRDLexample3}
\end{figure*}

Table~\ref{tab:PRDLresults} displays the CPU times, PSNR- and
SSIM-values and mean patch representation vector sparsity levels
obtained for the various measurement types, averaged over the instance
groups of the same size. The concrete examples in
Figures~\ref{fig:PRDLexample2} and \ref{fig:PRDLexample3} show the
results from DOLPHIn and plain Wirtinger Flow (WF; the real-valued,
$[0,1]$-box constrained variant, which corresponds to running
Algorithm~\ref{alg:PRDL} with $\mu=0$ and omitting the updates of $\A$
and $\D$). In all tests, we let the Wirtinger Flow method run for the
same number of iterations (75) and use the same starting points as for
the DOLPHIn runs.  Note that instead of random $\X_{(0)}$, we could
also use a spectral initialization similar to the one proposed for the
(unconstrained) Wirtinger Flow algorithm, see~\cite{CLS14}. Such
initialization can improve WF reconstruction (at least in the
noiseless case), and may also provide better initial estimates for
DOLPHIn. We have experimented with such a spectral approach and found
the results comparable to what is achievable with random $\X_{(0)}$,
both for WF and DOLPHIn. Therefore, we do not report these experiments
in the paper.

\begin{figure*}[tb]
\centering
\subfigure[]{\includegraphics[width=0.15\textwidth,resolution=300]{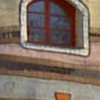}} 
$\qquad$
\subfigure[]{\includegraphics[width=0.15\textwidth,resolution=300]{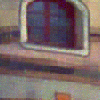}}
$\qquad$
\subfigure[]{\includegraphics[width=0.15\textwidth,resolution=300]{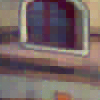}}
$\qquad$
\subfigure[]{\includegraphics[width=0.15\textwidth,resolution=300]{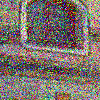}}
\caption{DOLPHIn example ``Waldspirale'' image, zoomed-in $100\times
  100$ pixel parts (magnified). (a) original image, (b) image
  reconstruction $\X_{\text{DOLPHIn}}$, (b) reconstruction $\cR(\D\A)$
  from sparsely coded patches, (c) reconstruction $\X_{\text{WF}}$
  after $75$ WF iterations. The slight block artefacts visible in (b)
  and (c) are due to the nonoverlapping patch approach in experiments
  and could easily be mitigated by introducing some patch overlap
  (cf., e.g., Fig.~\ref{fig:PRDLCexample}).}
\label{fig:PRDLexample3zoom}
\end{figure*} 

The DOLPHIn method consistently provides better image reconstructions
than WF, which clearly shows that our approach successfully introduces
sparsity into the phase retrieval problem and exploits it for
estimating the solution. As can be seen from
Table~\ref{tab:PRDLresults}, the obtained dictionaries allow for
rather sparse representation vectors, with the effect of making better
use of the information provided by the measurements, and also
denoising the image along the way. The latter fact can be seen in the
examples (Fig.~\ref{fig:PRDLexample2} and~\ref{fig:PRDLexample3}, see
also Fig.~\ref{fig:PRDLexample3zoom}) and also inferred from the
significantly higher PSNR and SSIM values for the estimates
$\X_{\text{DOLPHIn}}$ and $\cR(\D\A)$ (or $\proj_{\cX}(\cR(\D\A))$,
resp.) obtained from DOLPHIn compared to the reconstruction
$\X_{\text{WF}}$ of the WF algorithm (which cannot make use of hidden
sparsity). The gain in reconstruction quality is more visible in the
example of Fig.~\ref{fig:PRDLexample3}
(cf. Fig.~\ref{fig:PRDLexample3zoom}) than for that in
Fig.~\ref{fig:PRDLexample2}, though both cases assert higher
quantitative measures. Furthermore, note that DOLPHIn naturally has
higher running times than WF, since it performs more work per
iteration (also, different types of measurement operators require
different amounts of time to evaluate). Note also that storing $\A$
and $\D$ instead of an actual image $\X$ (such as the WF
reconstruction) requires saving only about half as many numbers
(including integer index pairs for the nonzero entries in $\A$).

As indicated earlier, the reconstruction $\cR(\D\A)$ is quite often
better than $\X_{\text{DOLPHIn}}$ w.r.t. at least one of either PSNR or SSIM
value. Nonetheless, $\X_{\text{DOLPHIn}}$ may be visually more appealing than
$\cR(\D\A)$ even if the latter exhibits a higher
quantitative quality measure (as is the case, for instance, in the
example of Figures~\ref{fig:PRDLexample3} and~\ref{fig:PRDLexample3zoom}); 
Furthermore, occasionally $\X_{\text{DOLPHIn}}$ achieves notably
better (quantitative) measures than $\cR(\D\A)$; an intuitive explanation 
may be that if, while the sparse coding of patches served
well to eliminate the noise and---by means of the patch-fit objective
term---to successfully ``push'' the $\X$-update steps toward a
solution of good quality, that solution eventually becomes ``so good'',
then the fact that $\cR(\D\A)$ is designed to be 
only an \emph{approximation} (of $\X$) predominates. 

On the other hand,
 $\X_{\text{DOLPHIn}}$ is sometimes
very close to $\X_{\text{WF}}$, which indicates a suboptimal setting
of the parameters $\mu$ and $\lambda$ that control how much
``feedback'' the patch-fitting objective term introduces into the
Wirtinger-Flow-like $\X$-update in the DOLPHIn algorithm. We discuss
parameter choices in more detail in the following subsection.

\begin{figure*}[tb]
\centering
\subfigure[]{\includegraphics[width=0.32\textwidth,resolution=300]{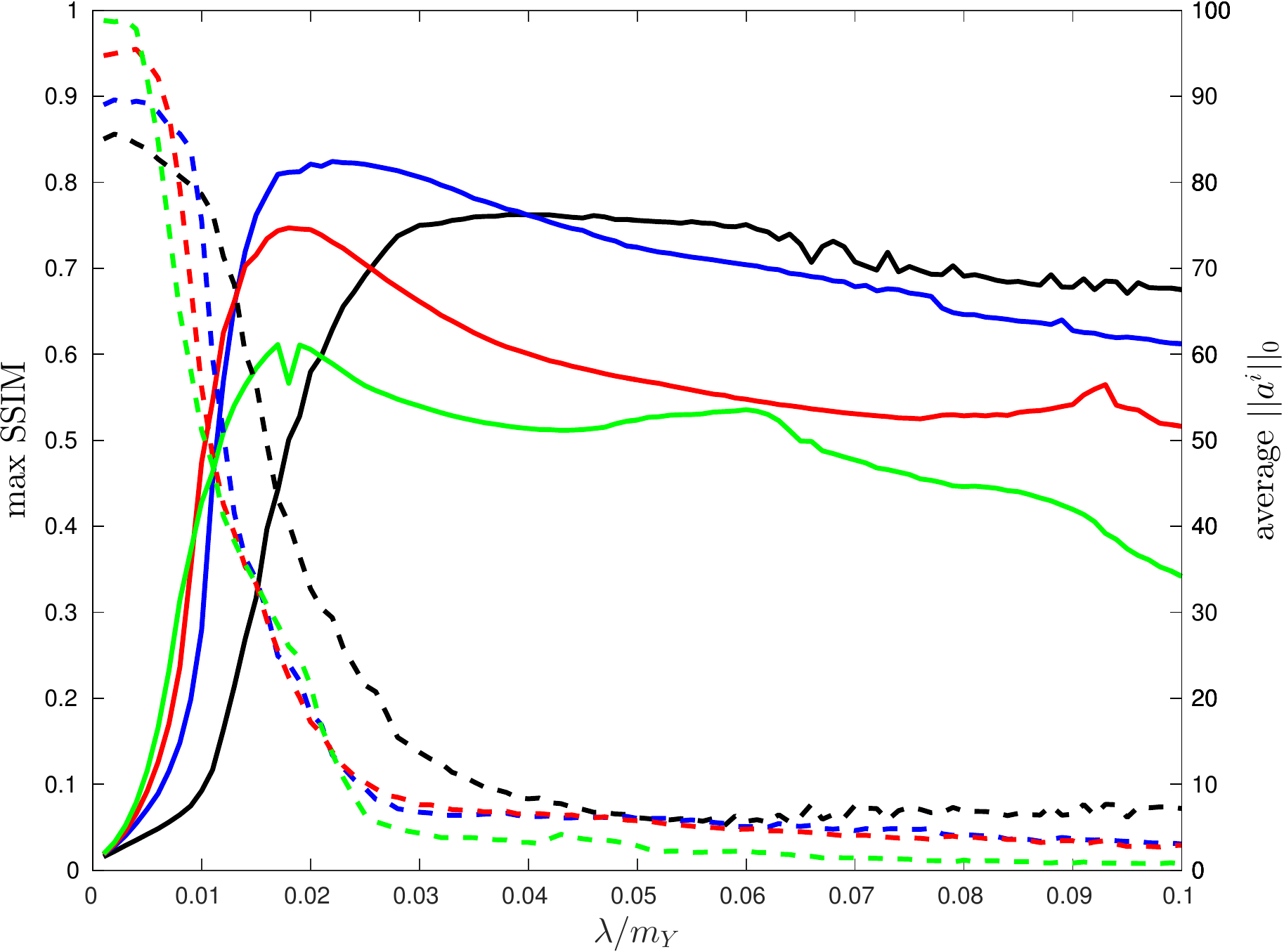}}
~
\subfigure[]{\includegraphics[width=0.32\textwidth,resolution=300]{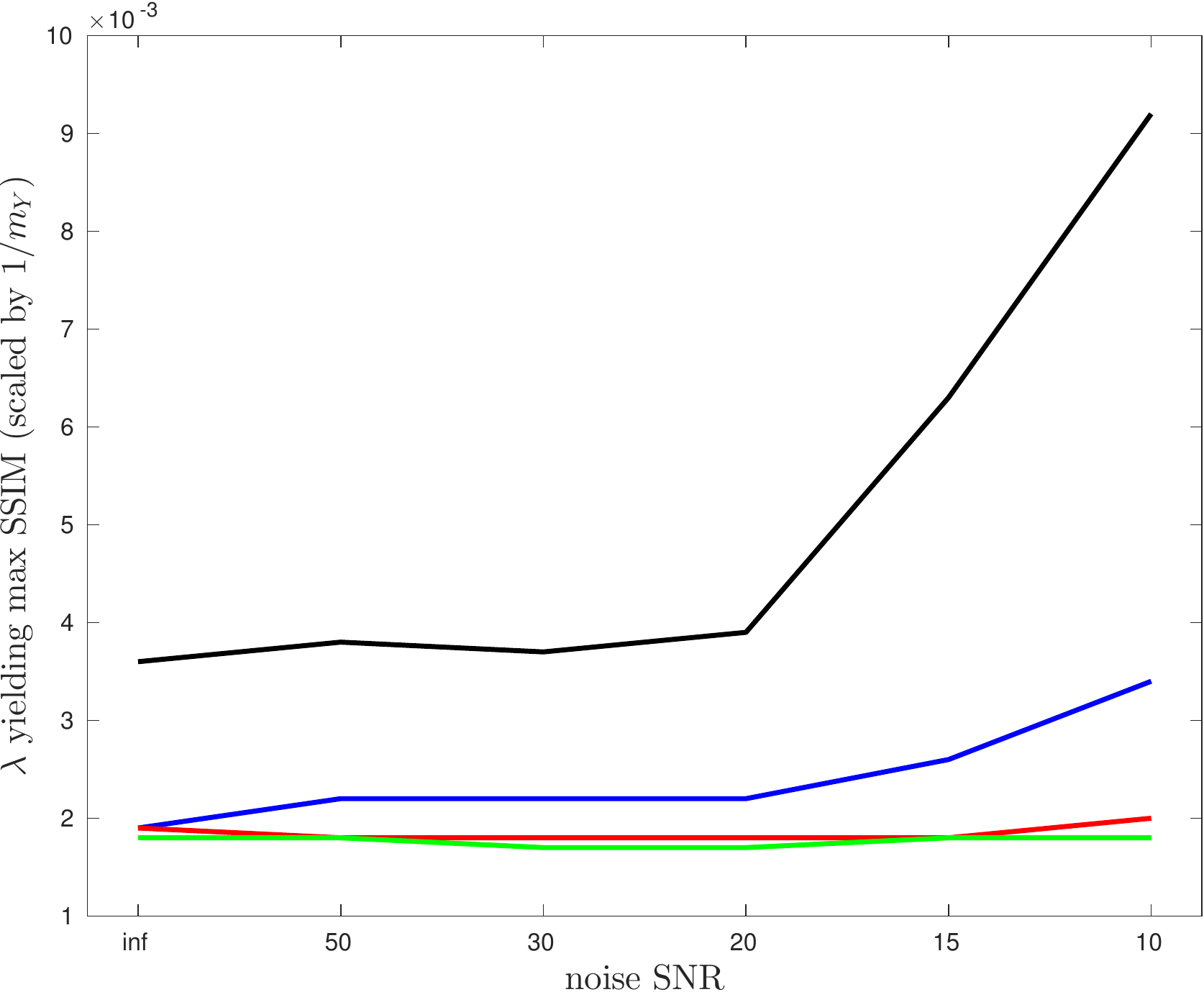}}
~
\subfigure[]{\includegraphics[width=0.32\textwidth,resolution=300]{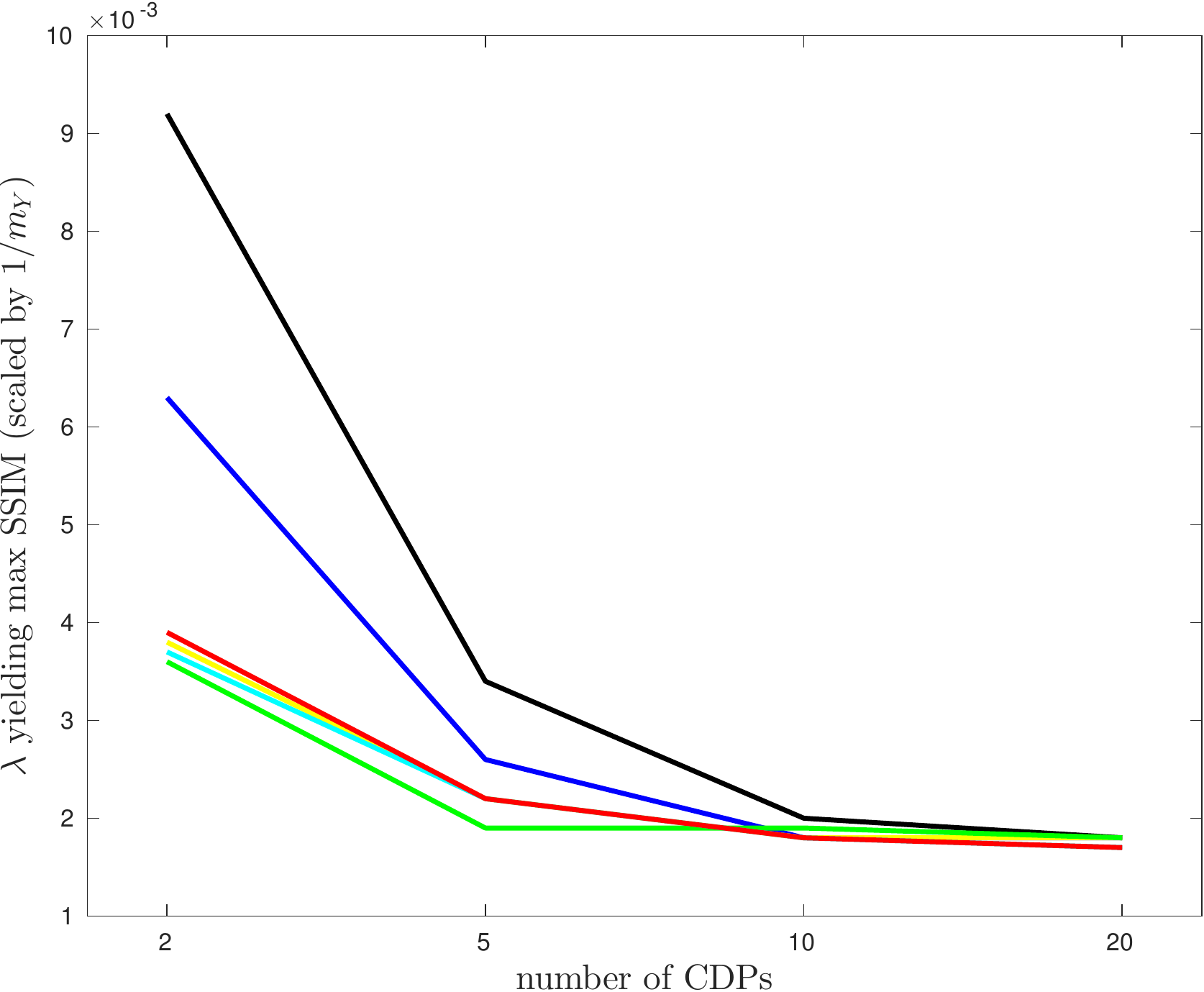}}
~
\subfigure[]{\includegraphics[width=0.32\textwidth,resolution=300]{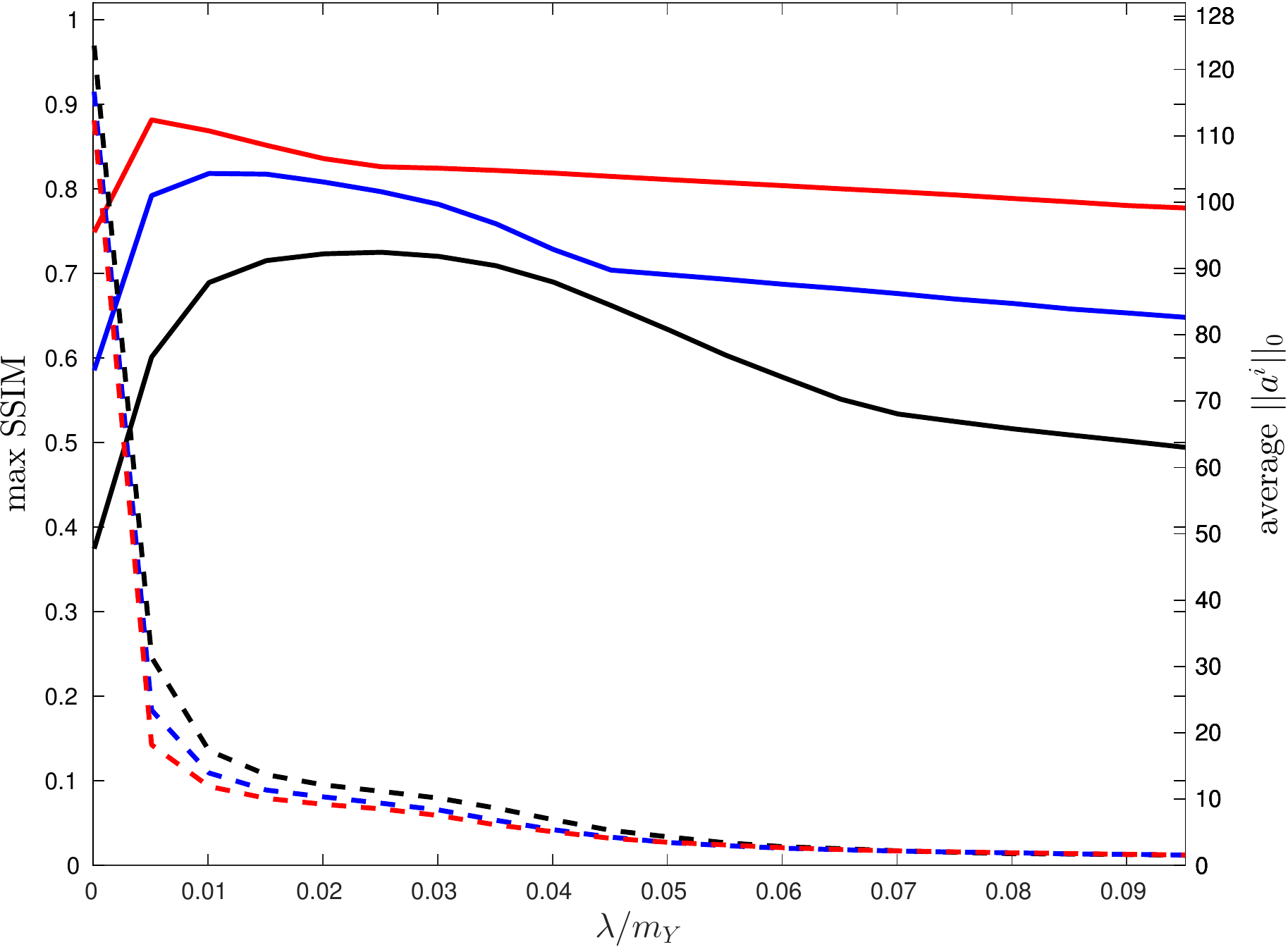}}
~
\subfigure[]{\includegraphics[width=0.32\textwidth,resolution=300]{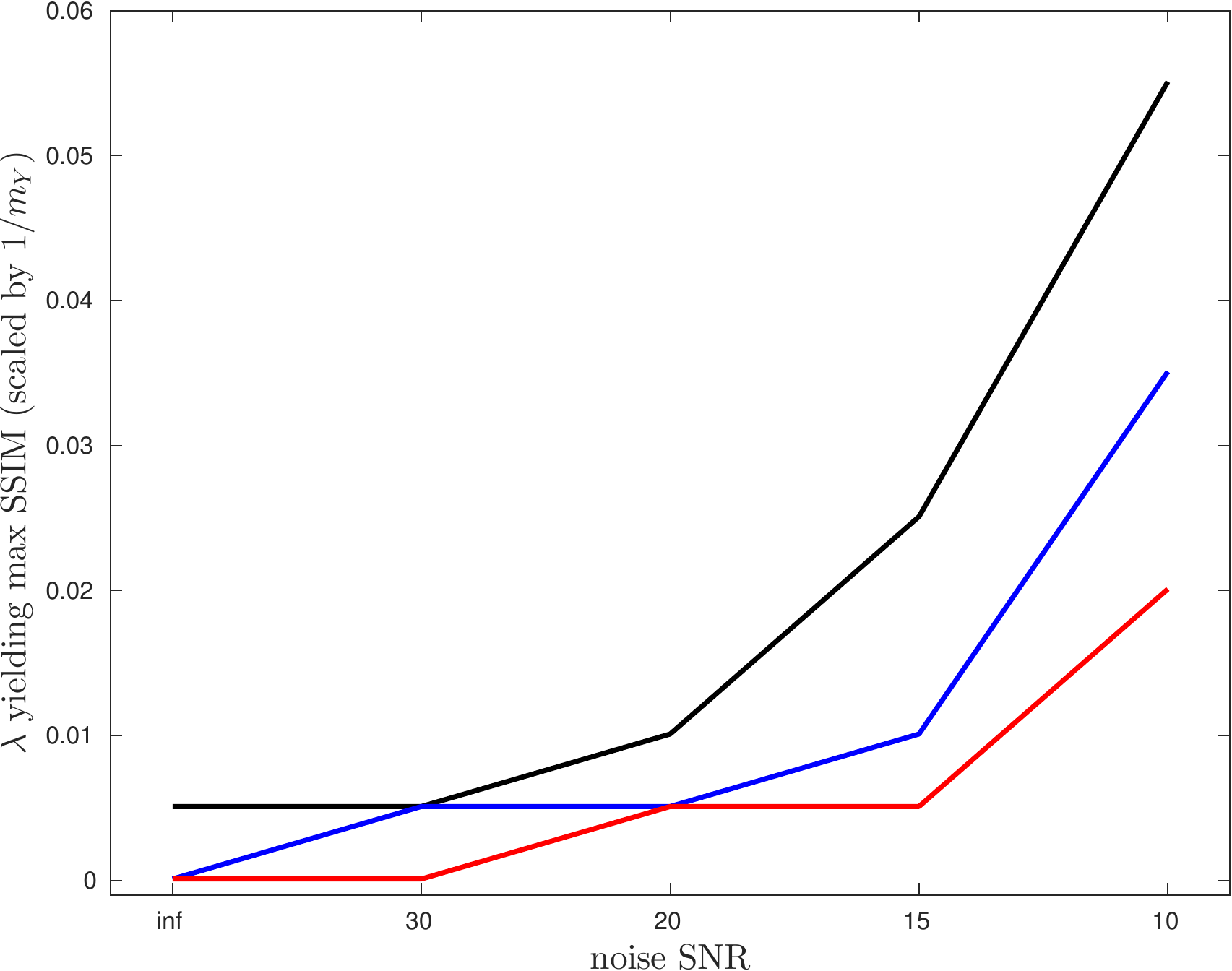}}
~
\subfigure[]{\includegraphics[width=0.32\textwidth,resolution=300]{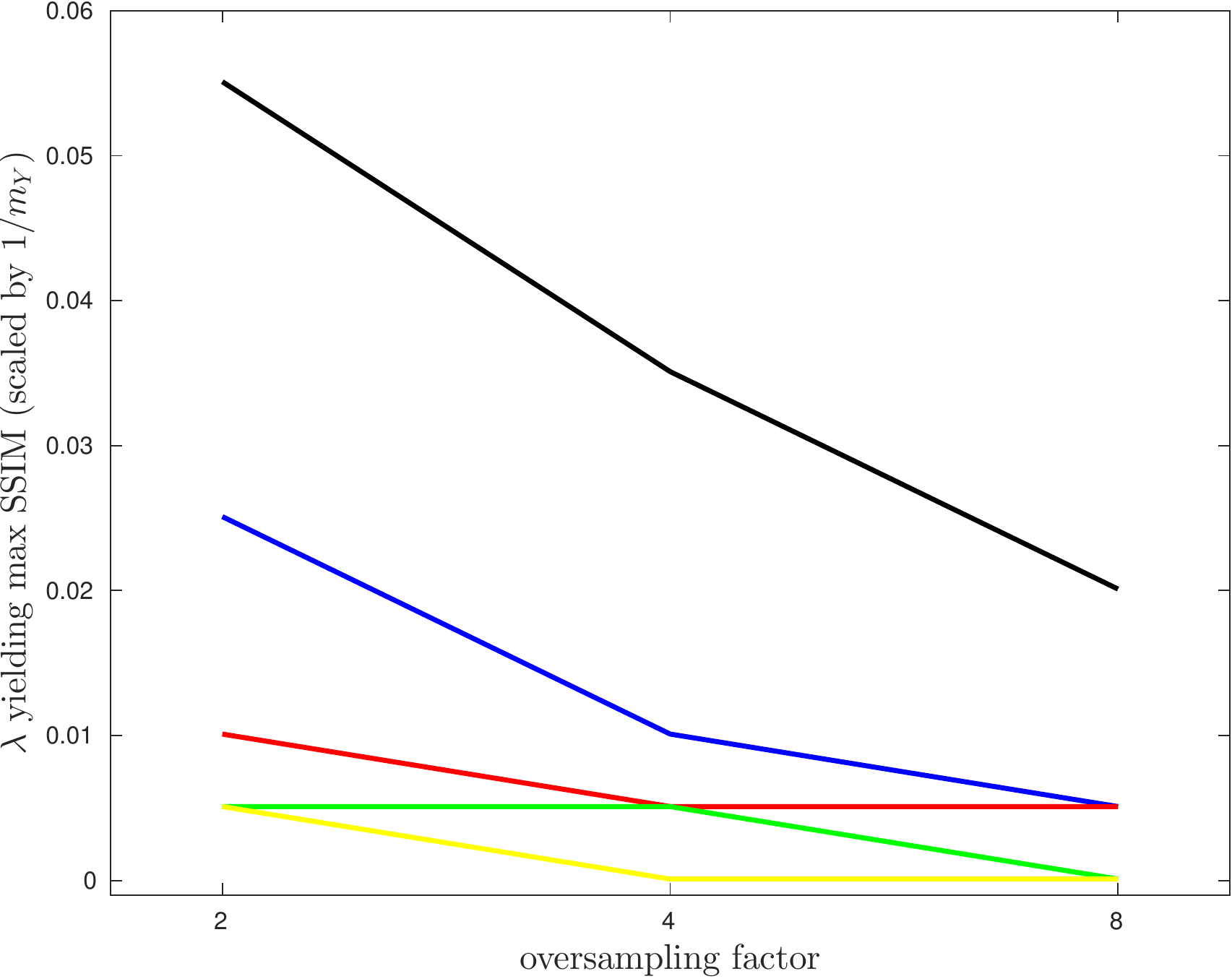}}

\caption{Influence of parameter $\lambda$ on best achievable SSIM
  values for reconstructed images and sensitivity w.r.t. sampling
  ratios and noise levels, for different measurement types. Fixed
  other parameters: $\mu=0.1 m_Y$, $K_1=25$, $K_2=50$, $s_1=s_2=8$
  (nonoverlapping patches), $\D_{(0)}=(\I,\F_D)$, $\X_{(0)}\in\cX$
  random. (a)--(c): Averages over reconstructions from ternary CDP
  measurements of the three $256\times 256$ images. The plots
  show (a) the best achievable SSIM for
  $\lambda/m_Y\in\{0.0001,0.0002,\dots,0.01\}$ (left vertical axis, solid
  lines) and average patch-sparsity of corresponding solutions (right
  vertical axis, dashed lines) for noise level
  SNR$(\Y,\abs{\cF(\hat{\X})}^2)=20$\,dB and (b) choice of $\lambda$
  yielding best SSIM for different noise levels, for number of masks
  $2$ (black), $5$ (blue), $10$ (red) or $20$ (green), respectively;
  (c) choice of $\lambda$ to achieve best SSIM for different number of
  masks, for noise-SNRs $10$ (black), $15$ (blue), $20$ (red), $30$
  (yellow), $50$ (light blue) and $\infty$ (green), respectively.
  (d)--(f): Averages over reconstructions from Gaussian measurements
  ($\Y=\abs{\G\hat{\X}}^2$) of the five $512\times 512$
  images. The plots display the same kind of information as (a)--(c),
  but in (d) with $\lambda/m_Y\in\{0.0001,0.0051,0.0101,\dots,0.0951\}$
  for noise-SNR $15$\,dB and in (e) for different noise levels, for
  sampling ratios $M_1/N_1=2$ (black), $4$ (blue) and $8$ (red),
  respectively; and in (f) with $M_1/N_1=2$ for noise-SNRs $10$
  (black), $15$ (blue), $20$ (red), $30$ (green) and $\infty$
  (yellow).}
\label{fig:ChoiceOfLambda}
\end{figure*}

\subsection{Hyperparameter Choices and Sensitivity}\label{subsec:parameters}
The DOLPHIn algorithm requires several parameters to be specified a
priori. Most can be referred to as \emph{design parameters}; the most
prominent ones are the size of image patches ($s_1\times s_2$),
whether patches should overlap or not (not given a name here), and the
number $n$ of dictionary atoms to learn. Furthermore, there are
certain \emph{algorithmic parameters} (in a broad sense) that need to
be fixed, e.g., the iteration limits $K_1$ and $K_2$ or the initial
dictionary $\D_{(0)}$ and image estimate $\X_{(0)}$. The arguably most
important parameters, however, are the model or \emph{regularization
  parameters} $\mu$ and $\lambda$. For any fixed combination of design
and algorithmic parameters in a certain measurement setup (fixed
measurement type/model and (assumed) noise level), it is conceivable
that one can find some values for $\mu$ and $\lambda$ that work well
for most instances, while the converse---choosing, say, iteration
limits for fixed $\mu$, $\lambda$ and other parameters---is clearly
not a very practical approach.

As is common for regularization parameters, a good ``general-purpose''
way to choose $\mu$ and $\lambda$ a priori is unfortunately not known.
To obtain the specific choices used in our experiments, we fixed all
the other parameters (including noise SNR and oversampling ratios),
then (for each measurement model) ran preliminary tests to identify
values $\mu$ for which good results could be produced with some
$\lambda$, and finally fixed $\mu$ at such values and ran extensive
benchmark tests to find $\lambda$ values that give the best results.

For DOLPHIn, $\mu$ offers some control over how much ``feedback'' from
the current sparse approximation of the current image estimate is
taken into account in the update step to produce the next image
iterate---overly large values hinder the progress made by the
Wirtinger-Flow-like part of the $\X$-update, while too small values
marginalize the influence of the approximation $\cR(\D\A)$, with one
consequence being that the automatic denoising feature is essentially
lost. Nevertheless, in our experiments we found that DOLPHIn is not
strongly sensitive to the specific choice of $\mu$ once a certain
regime has been identified in which one is able to produce meaningful
results (for some choice of $\lambda$). Hence, $\lambda$ may be
considered the most important parameter; note that this intuitively
makes sense, as $\lambda$ controls how strongly sparsity of the patch
representations is actually promoted, the exploitation of which to
obtain improved reconstructions being the driving purpose behind our
development of DOLPHIn.

Figure~\ref{fig:ChoiceOfLambda} illustrates the sensitivity of DOLPHIn
with respect to $\lambda$, in terms of reconstruction quality and
achievable patch-sparsity, for different noise levels, and examplary
measurement types and problem sizes. (In this figure, image quality is
measured by SSIM values alone; the plots using PSNR instead look very
similar and were thus omitted. Note, however, that the parameters
$\lambda$ yielding the best SSIM and PSNR values, respectively, need
not be the same.) As shown by (a) and (d), there is a clear
correlation between the best reconstruction quality that can be
achieved (in noisy settings) and the average sparsity of the patch
representation vectors $\a^i$. For larger noise, clearly a larger
$\lambda$ is needed to achieve good results---see (b) and (e)---which
shows that a stronger promotion of patch-sparsity is an adequate
response to increased noise, as is known for linear sparsity-based
denoising as well. Similarly, increasing the number of measurements
allows to pick a smaller $\lambda$ whatever the noise level actually
is, as can be seen in (c) and (f), respectively. The dependence of the
best $\lambda$ on the noise level appears to follow an exponential
curve (w.r.t. the reciprocal SNR) which is ``dampened'' by the
sampling ratio, i.e., becoming less steep and pronounced the more
measurements are available, cf. (b) and (e). Indeed, again referring
to the subplots (c) and (f), at a fixed noise level the best $\lambda$
values seem to decrease exponentially with growing number of
measurements.  It remains subject of future research to investigate
these dependencies in more detail, e.g., to come up with more or less
general (functional) rules for choosing~$\lambda$.

\subsection{Impact of Increased Inner Iteration Counts}\label{subsec:inneriters}
It is worth considering whether more inner iterations---i.e.,
consecutive update steps for the different variable blocks---lead to
further improvements of the results and\,/\,or faster convergence. In
general, this is an open question for block-coordinate descent algorithms, so the choices are typically made
empirically. Our default choices of $a=1$ ISTA iterations for the
$\A$-update (Step~\ref{step:PRDL_Aupdate} in
Algorithm~\ref{alg:PRDL}), $x=1$ projected gradient descent steps for
the $\X$-update (Step~\ref{step:PRDL_Xupdate}) and $d=1$ iterations of
the BCD scheme for the $\D$-update
(Steps~\ref{step:PRDL_DupdateStart}--\ref{step:PRDL_DupdateEnd})
primarily reflect the desire to keep the overall iteration
cost low. To assess whether another choice might yield significant
improvements, we evaluated the DOLPHIn performance for all
combinations of $a\in\{1,3,5\}$ and $d\in\{1,3,5\}$, keeping all other
parameters equal to the settings from the experiments reported on
above. (We also tried these combinations together with an increased
iteration count for the $\X$-update, but already for $x=2$ or $x=3$
the results were far worse than with just $1$ projected gradient
descent step; the reason can likely be found in the fact that without
adapting $\A$ and $\D$ to a modified $\X$-iterate, the patch-fitting
term of the objective tries to keep $\X$ close to a then-unsuitable
estimate $\cR(\D\A)$ based on older $\X$-iterates, which apparently
has a quite notable negative influence on the achievable progress in
the $\X$-update loop.) 

The results are summarized in condensed format
in Table~\ref{tab:InnerIterResults}, from which we can read off the
spread of the best and worst results (among the best ones achievable
with either $\X$ or $\proj_{\cX}(\cR(\D\A))$) for each
measurement-instance combination among all combinations
$(a,d)\in\{1,3,5\}^2$. (Full tables for each test run can be found
alongside our DOLPHIn code on the first author's webpage.) As the
table shows, the results are all quite close; while some
settings lead to sparser patch representations, the overall quality of
the best reconstructions for the various combinations usually differ
only slightly, and no particular combination stands out clearly as
being better than all others. In particular, comparing the results
with those in Table~\ref{tab:PRDLresults}, we find that our default
settings provide consistently good results; they may be improved upon
with some other combination of $a$ and $d$, but at least with the same
total iteration horizon ($K_1+K_2=75$), the improvements are often
only marginal. Since the overall runtime reductions (if any) obtained
with other choices for $(a,d)$ are also very small, there does not
seem to be a clear advantage to using more than a single iteration for
either update problem.

\begin{table*}[bt]
  \caption{Test results for DOLPHIn variants with different combinations of inner iteration numbers $a$ and $d$ for the $\A$- and $\D$-updates, resp. 
	Reported are the best mean values achievable (via either $\X_{\text{DOLPHIn}}$ or $\proj_{\cX}(\cR(\D\A))$) with any of the considered combinations $(a,d)\in\{1,3,5\}\times\{1,3,5\}$ (first rows for each measurement type) and the worst among the best values for each combination (second rows), along with the respective combinations yielding the stated values. All other parameters are identical to those 
	used in the experiments for default DOLPHIn ($a=d=1$), cf. Table~\ref{tab:PRDLresults}.
	}
  \label{tab:InnerIterResults}
  \centering
  \begin{tabular*}{\textwidth}{@{\extracolsep{\fill}}l@{\quad}c@{\quad}c@{\quad}c@{\quad}c@{\quad}c@{\quad}c@{\quad}c@{\quad}c@{\quad}c@{\quad}}
    \toprule
                                                 & \multicolumn{4}{c}{$256\times 256$ instances} & \multicolumn{4}{c}{$512\times 512$ instances}\\\cmidrule(r{9pt}){2-5}\cmidrule(r{9pt}){6-9}
    $\cF$ type & time & PSNR  & SSIM   & $\varnothing$\,$\norm{\a^i}_0$ & time  & PSNR  & SSIM   & $\varnothing$\,$\norm{\a^i}_0$\\
    \midrule
    $\G\hat\X$                 & 12.61 $(1,3)$ & 24.72 $(1,5)$ & 0.6680 $(1,5)$ & 3.69 $(1,5)$ & 68.18 $(1,3)$ & 24.43 $(1,3)$ & 0.6903 $(3,3)$ & 6.25 $(1,5)$\\
                               & 15.30 $(5,5)$ & 24.48 $(3,3)$ & 0.6485 $(3,3)$ & 4.94 $(3,1)$ & 79.28 $(5,5)$ & 24.30 $(5,3)$ & 0.6803 $(1,5)$ & 8.55 $(3,5)$\\[0.5em]
    $\G\hat\X\G^*$             & 51.02 $(1,3)$ & 23.71 $(1,5)$ & 0.7330 $(1,5)$ & 6.71 $(5,5)$ & 357.49 $(1,1)$ & 23.44 $(1,5)$ & 0.7685 $(1,1)$ & 9.54 $(5,5)$\\
                               & 54.98 $(3,1)$ & 23.57 $(5,5)$ & 0.7188 $(1,3)$ & 7.55 $(1,3)$ & 371.20 $(5,3)$ & 23.39 $(5,1)$ & 0.7633 $(1,3)$ & 11.61 $(1,3)$\\[0.5em]
    $\G\hat\X\H^*$             & 50.95 $(1,3)$ & 23.68 $(1,1)$ & 0.7315 $(1,1)$ & 6.64 $(5,\{1,5\})$ & 357.66 $(1,1)$ & 23.44 $(1,\{3,5\})$ & 0.7693 $(1,3)$ & 9.46 $(5,5)$\\
                               & 56.05 $(3,1)$ & 23.56 $(5,3)$ & 0.7143 $(5,3)$ & 7.59 $(1,3)$ & 373.43 $(5,5)$ & 23.40 $(5,\{3,5\})$ & 0.7650 $(5,3)$ & 11.57 $(1,5)$\\[0.5em]
    CDP (cf.~\eqref{eq:CDPop}) & 8.56 $(1,1)$ & 27.19 $(3,1)$ & 0.7692 $(3,1)$ & 7.71 $(1,3)$ & 35.66 $(1,3)$ & 27.38 $(5,1)$ & 0.7837 $(3,3)$ & 11.40 $(1,3)$ \\
                               & 11.40 $(5,5)$ & 27.04 $(3,5)$ & 0.7654 $(1,1)$ & 9.49 $(3,1)$ & 47.07 $(5,3)$ & 27.33 $(1,\{1,3\})$ & 0.7818 $(1,5)$ & 13.43 $(3,1)$\\
  \bottomrule
  \end{tabular*}
\end{table*}

\subsection{Influence of the First DOLPHIn Phase}\label{subsec:firstphase}
Our algorithm keeps the dictionary fixed at its initialization for the
first $K_1$ iterations in order to prevent the dictionary from
training on relatively useless first reconstruction iterates. Indeed,
if all variables \emph{including $\D$} are updated right from the
beginning (i.e., $K_1=0$, $K_2=75$), then we end up with inferior results
(keeping all other parameters unchanged): The obtained patch
representations are much less sparse, the quality of the image
estimate $\cR(\D\A)$ decreases drastically, and also the
quality of the reconstruction $\X$ becomes notably worse. This
demonstrates that the dictionary apparently ``learns too much noise''
when updated from the beginning, and the positive effect of filtering
out quite a lot of noise in the first phase by regularizing with
sparsely coded patches using a fixed initial dictionary is almost
completely lost. To give an example, for the CDP testset on $256\times
256$ images, the average values obtained by DOLPHIn when updating also
the dictionary from the first iteration onward are: $9.24$ seconds
runtime (versus $8.56$ for default DOLPHIn,
cf. Table~\ref{tab:PRDLresults}), mean patch sparsity
$\varnothing\norm{\a^i}_0\approx 20.09$ (vs. $7.85$), PSNR
$26.80$\,dB and $8.88$\,dB (vs. $27.15$ and $26.58$) and SSIM-values
$0.6931$ and $0.0098$ (vs. $0.7416$ and $0.7654$) for the
reconstructions $\X$ and $\proj_{\cX}(\cR(\D\A))$, respectively. 

On the other hand, one might argue that if the first iterates are
relatively worthless, the effort of updating $\A$ in the first DOLPHIn
phase (i.e., the first $K_1$ iterations) could be saved as
well. However, the influence of the objective terms involving $\D$ and
$\A$ should then be completely removed from the algorithm for the
first $K_1$ iterations; otherwise, the patch-fitting term will
certainly hinder progress made by the $\X$-update because it then
amounts to trying to keep $\X$ close to the initial estimate $\cR(\D\A)$,
which obviously needs not bear any resemblance to the sought solution
at all. Thus, if both $\D$ and $\A$ are to be unused in the first
phase, then one should temporarily set $\mu=0$, with the consequence that
the first phase reduces to pure projected gradient descent for $\X$
with respect to the phase retrieval objective---i.e., essentially, Wirtinger
Flow. Therefore, proceeding like this simply amounts to a different
initialization of $\X$. Experiments with this DOLPHIn variant
($K_1=25$ initial WF iterations followed by $K_2=50$ full DOLPHIn
iterations including $\A$ and $\D$; all other parameters again left
unchanged) showed that the achievable patch-sparsities remain about
the same for the Gaussian measurement types but become much worse for
the CDP setup, and that the average PSNR and SSIM values become (often
clearly) worse in virtually all cases. In the example of measurements
$\G\hat{\X}$ of the $512\times 512$ test images, the above-described
DOLPHIn variant runs $62.34$ seconds on average (as less effort is
spent in the first phase, this is naturally lower than the $68.62$
seconds default DOLPHIn takes), produces slightly lower average
patch-sparsity ($5.88$ vs. $6.30$ for default DOLPHIn), but for both
$\X$ and $\proj_{\cX}(\cR(\D\A))$, the PSNR and SSIM values are
notably worse ($22.32$\,dB and $20.55$\,dB vs. $24.42$\,dB and
$22.66$\,dB, and $0.6165$ and $0.5921$ vs. $0.6547$ and $0.6807$,
resp.). The reason for the observed behavior can be found in the
inferior performance of WF without exploiting patch-sparsity (cf. also
Table~\ref{tab:PRDLresults}); note that the results also further
demonstrate DOLPHIn's robustness w.r.t. the initial
point---apparently, the initial point obtained from some WF iterations
is not more helpful for DOLPHIn than a random first guess.

Finally, it is also worth considering what happens if the dictionary
updates are turned off completely, i.e., $K_2=0$. Then, DOLPHIn
reduces to patch-sparsity regularized Wirtinger Flow, a WF variant
that apparently was not considered previously. Additional experiments
with $K_1=75$, $K_2=0$ and $\D=\D_{(0)}=(\I,\F_D)$ (other parameters
left unchanged) showed that this variant consistently produces higher
sparsity (i.e., smaller average number of nonzero entries) of the
patch representation coefficient vectors, but that the best
reconstruction ($\X$ or $\proj_{\cX}(\cR(\D\A))$) is always
significantly inferior to the best one produced by our default version
of DOLPHIn. The first observation is explained by the fact that with
$\D$ fixed throughout, the patch coding ($\A$-update) only needs to
adapt w.r.t. new $\X$-iterates but not a new $\D$; at least if the new
$\X$ is not too different from the previous one, the former
representation coefficient vectors still yield an acceptable
approximation, which no longer holds true if the dictionary was also
modified. While it should also me mentioned that the patch-sparsity
regularized version performs much better than plain WF already, the
second observation clearly indicates the additional benefit of working
with trained dictionaries, i.e., superiority of DOLPHIn also over the
regularized WF variant.

Full tables containing results for all testruns reported on in
this subsection are again available online along with our DOLPHIn
code on the first author's website.

\subsection{Sparsity-Constrained DOLPHIn Variant}\label{subsec:l0prdl}
\begin{table*}[tb]
  \caption{Test results for sparsity-constrained DOLPHIn, using overlapping patches, for $m_Y$ Gaussian-type and coded diffraction pattern (CDP) measurements. Reported are mean values (geometric mean for CPU times) per measurement type, obtained from three instances with random $\X_{(0)}$ and random noise for each of the three $256\times 256$ and five $512\times 512$ images, w.r.t. the reconstructions from DOLPHIn ($\X_{\text{DOLPHIn}}$ and $\proj_{\cX}(\cR(\D\A))$) with parameters $(\mu_1,\mu_2)$ and (real-valued, $[0,1]$-constrained) Wirtinger Flow ($\X_{\text{WF}}$), respectively. (CPU times in seconds, PSNR in decibels).}
  \label{tab:PRDLCresults}
  \centering
\begin{tabular*}{\textwidth}{@{\extracolsep{\fill}}l@{\quad}l@{\quad}c@{\quad}c@{\quad}c@{\quad}c@{\quad}c@{\quad}c@{\quad}c@{\quad}c@{\quad}c@{\quad}c@{\quad}c@{\quad}}
\toprule
 & & \multicolumn{5}{c}{$256\times 256$ instances} & \multicolumn{5}{c}{$512\times 512$ instances}\\\cmidrule(r{9pt}){3-7}\cmidrule(r{9pt}){8-12}
$\cF$ type & reconstruction & $(\mu_1,\mu_2)/m_Y$ & time & PSNR & SSIM & $\varnothing$\,$\norm{\a^i}_0$ & $(\mu_1,\mu_2)/m_Y$ & time & PSNR & SSIM & $\varnothing$\,$\norm{\a^i}_0$\\
\midrule
$\G\hat\X$ & $\X_{\text{DOLPHIn}}$ & (0.005,0.0084) & 23.58 & 26.79 & 0.6245 & & (0.005,0.0084) & 137.15 & 26.73 & 0.7090 & \\
& $\proj_{\cX}(\cR(\D\A))$ & & & 27.98 & 0.7721 & 2.71 & & & 27.60 & 0.7786 & 3.28\\
& $\X_{\text{WF}}$ & & 5.04 & 19.00 & 0.2889 & -- & & 32.63 & 18.94 & 0.3811 & --\\[0.5em]
$\G\hat\X\G^*$ & $\X_{\text{DOLPHIn}}$ & (0.005,0.0084) & 57.52 & 22.71 & 0.4143 & & (0.005,0.0084) & 358.29 & 22.49 & 0.5234 & \\
& $\proj_{\cX}(\cR(\D\A))$ & & & 27.29 & 0.6129 & 7.96 & & & 27.47 & 0.7202 & 8.00\\
& $\X_{\text{WF}}$ & & 32.44 & 22.71 & 0.4145 & -- & & 232.26 & 22.49 & 0.5239 & --\\[0.5em]
$\G\hat\X\H^*$ & $\X_{\text{DOLPHIn}}$ & (0.005,0.0084) & 57.67 & 22.54 & 0.4082 & & (0.005,0.0084) & 356.27 & 22.56 & 0.5272 & \\
& $\proj_{\cX}(\cR(\D\A))$ & & & 27.14 & 0.6059 & 7.97 & & & 27.55 & 0.7233 & 8.00\\
& $\X_{\text{WF}}$ & & 32.32 & 22.56 & 0.4088 & -- & & 232.62 & 22.56 & 0.5276 & --\\[0.5em]
CDP (cf.~\eqref{eq:CDPop}) & $\X_{\text{DOLPHIn}}$ & (0.005,0.0084) & 22.10 & 28.00 & 0.8041 & & (0.005,0.0084) & 112.33 & 26.30 & 0.7400 & \\
& $\proj_{\cX}(\cR(\D\A))$ & & & 26.87 & 0.7789 & 1.52 & & & 24.95 & 0.6557 & 1.51\\
& $\X_{\text{WF}}$ & & 2.70 & 21.95 & 0.3880 & -- & & 13.68 & 21.73 & 0.4935 & --\\
\bottomrule
\end{tabular*}
\end{table*}

From Table~\ref{tab:PRDLresults} and Figure~\ref{fig:ChoiceOfLambda},
(a) and (d), it becomes apparent that a sparsity level of $8\pm 4$
accompanies the good reconstructions by DOLPHIn. This suggests use in
a DOLPHIn variant we already briefly hinted at: Instead of using
$\ell_1$-norm regularization, we could incorporate explicit sparsity
constraints on the $\a^i$. The corresponding DOLPHIn model then reads
\begin{align}
  \nonumber\min_{\X,\D,\A}\,&\tfrac14\Norm{\Y-\abs{\cF(\X)}^2}_{\text{F}}^2+\tfrac{\mu}{2}\bignorm{\cE(\X)-\D\A}_{\text{F}}^2\\
  \label{prob:PRDLC}\suchthat~~\,&\X\in[0,1]^{N_1\times N_2},~~\D\in\cD,~~\norm{\a^i}_0\leq k~\forall\,i=1,\dots,p,
\end{align}
where $k$ is the target sparsity level. We no longer need to tune the
$\lambda$ parameter, and can let our previous experimental results
guide the choice of $k$. Note that the only modification to
Algorithm~\ref{alg:PRDL} concerns the update of~$\A$
(Step~\ref{step:PRDL_Aupdate}), which now requires solving or
approximating
$$
\min\,\tfrac{1}{2}\bignorm{\cE(\X)-\D\A}_{\text{F}}^2\quad\suchthat\quad\norm{\a^i}_0\leq k\quad\forall\,i=1,\dots,p.\\
$$
In our implementation, we do so by running Orthogonal Matching Pursuit
(OMP) \cite{PRK93} for each column $\a^i$ of $\A$ separately until
either the sparsity bound is reached or
$\norm{\x^i-\D\a^i}_2\leq\varepsilon$, where we set a default of
$\varepsilon\define 0.1$. (The value of $\varepsilon$ is again a
parameter that might need thorough benchmarking; obviously, it is
related to the amount of noise one wants to filter out by sparsely
representing the patches---higher noise levels will require larger
$\varepsilon$ values. The $0.1$ default worked quite well in our
experiments, but could probably be improved by benchmarking as well.)
The OMP code we used is also part of the SPAMS package.

The effect of the parameter $\mu$ is more pronounced in the
sparsity-constrained DOLPHIn than in Algorithm~\ref{alg:PRDL}; however,
it appears its choice is less dependent on the measurement model
used. By just a few experimental runs, we found that good results in
all our test cases can be achieved using $K_1=K_2=25$ iterations,
where in the first $K_1$ (with fixed dictionary), we use a value of
$\mu=\mu_1=0.005 m_Y$ along with a sparsity bound $k=k_1=4$, and in
the final $K_2$ iterations (in which the dictionary is updated),
$\mu=\mu_2=1.68\mu_1=0.0084 m_Y$ and $k=k_2=8$. Results on the same
instances considered before (using the same $\D_{(0)}$) are presented
in Table~\ref{tab:PRDLCresults}, with the exception that for the CDP
case, we used $2$ complex-valued octanary masks (cf.~\cite{CLS14})
here; the initial image estimates and measurement noise were again
chosen randomly. Note also that for these test, we used complete sets
of \emph{overlapping} patches. This greatly increases $p$ and hence
the number of subproblems to be solved in the $\A$-update step, which
is the main reason for the increased running times compared to
Table~\ref{tab:PRDLresults} for the standard DOLPHIn method. (It should
however also be mentioned that OMP requires up to $k$ iterations per
subproblem, while in Algorithm~\ref{alg:PRDL} we explicitly restricted
the number of ISTA iterations in the $\A$-update to just a single one.)

\begin{figure*}[tb]
\centering
\subfigure[]{\includegraphics[width=0.23\textwidth,height=0.23\textwidth,resolution=300]{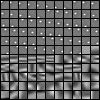}}
~
\subfigure[]{\includegraphics[width=0.23\textwidth,height=0.23\textwidth,resolution=300]{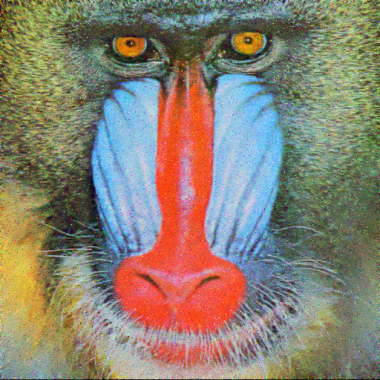}} 
~
\subfigure[]{\includegraphics[width=0.23\textwidth,resolution=300]{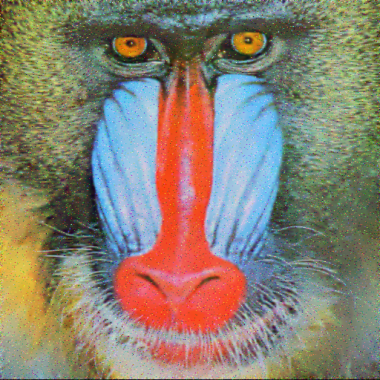}}
~
\subfigure[]{\includegraphics[width=0.23\textwidth,resolution=300]{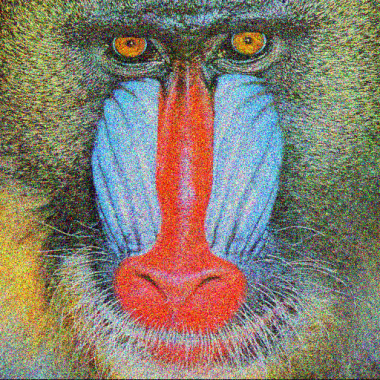}}\\
\subfigure[]{\includegraphics[width=0.15\textwidth,resolution=300]{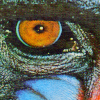}}
$\qquad$
\subfigure[]{\includegraphics[width=0.15\textwidth,resolution=300]{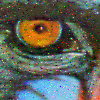}}
$\qquad$
\subfigure[]{\includegraphics[width=0.15\textwidth,resolution=300]{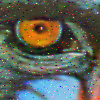}}
$\qquad$
\subfigure[]{\includegraphics[width=0.15\textwidth,resolution=300]{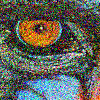}}

\caption{Sparsity-Constrained DOLPHIn example: Image original is the
  $512\times 512$ RGB ``mandrill'' picture, measurements are noisy
  CDPs (obtained using two complex-valued octanary masks, noise-SNR
  $10$\,dB, per color channel), $\mu_1=0.003\,m_Y$,
  $\mu_2=0.00504\,m_Y$ (other parameters as described in
  Section~\ref{subsec:l0prdl}). $\D_{(0)}=(\I,\F_D)$ for R-channel,
  final dictionary then served as initial dictionary for G-channel,
  whose final dictionary in turn was initial dictionary for B-channel;
  $\X_{(0)}\in\cX$ random for each channel. (a) final dictionary
  (excerpt) for B-channel (b) image reconstruction
  $\X_{\text{DOLPHIn}}$, (c) image reconstruction $\cR(\D\A)$ from
  sparsely coded patches, (d) reconstruction $\X_{\text{WF}}$ after
  $50$ WF
  iterations. 
  Final PSNR values: $20.53$\,dB for $\cR(\D\A)$, $20.79$\,dB for
  $\X_{\text{DOLPHIn}}$, $14.47$\,dB for $\X_{\text{WF}}$; final SSIM
  values: $0.4780$ for $\cR(\D\A)$, $0.5242$ for
  $\X_{\text{DOLPHIn}}$, $0.2961$ for $\X_{\text{WF}}$; average
  $\norm{\a^i}_0$ is $5.30$. (Means over all color channels.)
  (e)--(h): zoomed-in $100\times 100$ pixel parts (magnified) of (e)
  original image, (f) $\X_{\text{DOLPHIn}}$, (g) $\cR(\D\A)$ and (h)
  $\X_{\text{WF}}$.}
\label{fig:PRDLCexample}
\end{figure*}

A concrete example is given in Figure~\ref{fig:PRDLCexample}; here, we
consider the color ``mandrill'' image, for which the reconstruction
algorithms (given just two quite heavily noise-corrupted octanary CDP
measurements) were run on each of the three RGB channels
separately. The sparsity-constrained DOLPHIn reconstructions appear
superior to the plain WF solution both visually and in terms of the
quality measures PSNR and SSIM.

\section{Discussion and Conclusion}\label{sec:conclusion}
In this paper, we introduced a new method, called DOLPHIn, for
dictionary learning from noisy nonlinear measurements without phase
information. In the context of image reconstruction, the algorithm
fuses a variant of the recent Wirtinger Flow method for phase
retrieval with a patch-based dictionary learning model to obtain
sparse representations of image patches, and yields monotonic
objective decrease or (with appropriate step size selection)
convergence to a stationary point for the nonconvex combined DOLPHIn
model.

Our experiments demonstrate that dictionary learning
for phase retrieval with a patch-based sparsity is a promising
direction, especially for cases in which the original Wirtinger Flow
approach fails (due to high noise levels and/or relatively low
sampling rates).

Several aspects remain open for future research. For instance,
regarding the generally difficult task of parameter tuning, additional
benchmarking for to-be-identified instance settings of special
interest could give further insights on how to choose, e.g., the
regularization parameters in relation to varying noise levels.

It may also be worth developing further variants of our algorithm; the
successful use of $\ell_0$-constraints instead of
the~$\ell_1$-penalty, combining OMP with our framework, is just one
example. Perhaps most importantly, future research will be directed
towards the ``classic'' phase retrieval problem in which one is given
the (squared) magnitudes of Fourier measurements, see, e.g.,
\cite{SECCMS15,M07,F82}. Here, the WF 
method fails, and existing other (projection-based) methods are not
always reliable either. The hope is that introducing sparsity via a
learned dictionary will also enable improved reconstructions in the
Fourier setting.

To evaluate the quality of the learned dictionary, one might also ask
how DOLPHIn compares to the straightforward approach to first run
(standard) phase retrieval and then learn dictionary and sparse patch
representations from the result. Some preliminary experiments (see
also those in Section~\ref{subsec:firstphase} pertaining to keeping both
$\A$ and $\D$ fixed in the first DOLPHIn phase)
indicate 
that both approaches produce comparable results in the noisefree
setting, while our numerical results demonstrate a denoising feature
of our algorithm that the simple approach obviously lacks.

Similary, it will be of interest to see if the dictionaries learned by
DOLPHIn can be used successfully within sparsity-aware methods (e.g.,
the Thresholded WF proposed in~\cite{CLM15}, if that were modified to
handle local (patch-based) sparsity instead of global priors). In
particular, learning dictionaries for patch representations of images
from a whole \emph{class} of images would then be an interesting point
to consider. To that end, note that the DOLPHIn model and algorithm
can easily be extended to multiple input images whose patches
are all to be represented using a single dictionary by summing up
the objectives for each separate image, but with the same $\D$
throughout.

Another interesting aspect to evaluate is by how much                     
reconstruction quality and achievable sparsity degrade due to the loss
of phase (or, more generally, measurement nonlinearity), 
compared to the linear measurement case.

\appendices
\section{One-Dimensional DOLPHIn}\label{sec:appendixA}
In the following, we derive the DOLPHIn algorithm for the
one-dimensional setting~\eqref{eq:PR}. In particular, the (gradient)
formulas for the 2D-case can be obtained by applying the ones given
below to the vectorized image $\x=\vect(\X)$ (stacking columns of $\X$
on top of each other to form the vector $\x$), the vectorized matrix
$\a=\vect(\A)$, and interpreting the matrix $\F\in\C^{M\times N}$ as
describing the linear transformation corresponding to $\cF$ in terms
of the vectorized variables.

We now have a patch-extraction matrix $\P_e\in\R^{ps\times N}$ which
gives us $\P_e\x=((\x^1)^\top,\dots,(\x^p)^\top)^\top$ (in the
vectorized 2D-case, $\x^i$ then is the vectorized $i$-th patch of
$\X$, i.e., $\P_e$ corresponds to $\cE$). Similarly, we have a
patch-reassembly matrix $\P_a\in\R^{N\times ps}$; then, the reassembled
signal vector will be $\P_a((\a^1)^\top \D^\top,\dots,(\a^p)^\top
\D^\top)^\top$ (so $\P_a$ corresponds to $\cR$). Note that
$\P_a=\P_e^\dagger =\big(\P_e^\top\P_e\big)^{-1}\P_e^\top$; in
particular, $\x=\P_a\P_e\x$, and $\P_e^\top \P_e$ is a diagonal matrix
for which each diagonal entry is associated to a specific vector
component and gives the number of patches this component is contained
in. (Thus, if $\x=\vect(\X)$ is a vectorized 2D-image,
$\P_e^\top\P_e\x=\vect(\mR\odot\X)$ with $\mR$ as defined in
Section~\ref{subsec:Xupdate}.) Note that for nonoverlapping patches,
$\P_e$ is simply a permutation matrix, and $\P_a=\P_e^\top$ (so
$\P_a\P_e=\I$). Also, applying just $\P_e^\top$ actually reassembles a
signal from patches by simply adding the component's values
\emph{without} averaging.

We wish to represent each patch as $\x^i\approx\D\a^i$ with sparse
coefficient vectors $\a^i$; with
$\a\define((\a^1)^\top,\dots,(\a^p)^\top)^\top$ and
$\hat\D\define\I_p\otimes\D$, this sparse-approximation relation can
be written as $\P_e\x\approx\hat\D\a$. Our model to tackle the
1D-problem~\eqref{eq:PR} reads
\begin{align}
  \nonumber\min_{\x,\D,\al}\quad &\tfrac14\bignorm{\y-\abs{\F\x}^2}_2^2+ \tfrac{\mu}{2}\bignorm{\P_e\x-\hat\D\a}_2^2+\lambda\norm{\a}_1\\
  \label{prob:PRDL1d}\text{s.t.}\quad &\x\in\cX\define[0,1]^N,\quad\D\in\cD;
\end{align}
here, $\y\define\abs{\F\hat\x}^2+\n$, with $\hat\x$ the original signal
we wish to recover and $\n$ a noise vector.

The update formulas for $\a$ (separately for $\a^1,\dots,\a^p$) and
$\D$ remain the same as described before, see
Sections~\ref{subsec:Aupdate} and~\ref{subsec:Dupdate}, respectively.
However, the update problem for the phase retrieval solution---now
derived from~\eqref{prob:PRDL1d}, with $\D$ and $\a$ fixed at their
current values---becomes decreasing the objective
\begin{equation}\label{prob:PRDL1dXupdate}
  \tfrac14\bignorm{\y-\abs{\F\x}^2}_2^2+\tfrac{\mu}{2}\bignorm{\P_e\x-\hat\D\a}_2^2\quad\text{with}\quad\x\in\cX.
\end{equation}
We approach this by means of a projected gradient descent step; since
we consider real-valued $\x$-variables, this essentially amounts to
(one iteration of) the Wirtinger Flow method~\cite{CLS14},
accommodating the $[0,1]$-box constraints via projection onto them
after the (Wirtinger) gradient step. The step size will be chosen to
achieve a reduction of the objective w.r.t. its value at the previous
$\x$.

The objective function~\eqref{prob:PRDL1dXupdate} can be rewritten as
\[
\xi(\x) \define\phi(\x)+\psi(\x) \define\tfrac{1}{4}\sum_{j=1}^N\big(y_j-\x^\top\F_{j\cdot}^*\F_{j\cdot}\x\big)^2+\tfrac{\mu}{2}\norm{\P_{e}\x-\hat\D\a}_2^2.
\]
The gradient of $\psi(\x)$ is straightforwardly found to be
\[
   \nabla\psi(\x)=\mu\P_e^\top\big(\P_e\x-\hat\D\a\big).
\]
Regarding
$\nabla\phi(\x)=\tfrac{1}{2}\sum_{j=1}^N\big(y_j-\x^\top\M^j\x\big)\cdot\nabla
\big(y_j-\x^\top \M^j\x\big)$, where $\M^j\define
\F_{j\cdot}^*\F_{j\cdot}$ ($\F_{j\cdot}$ is the $j$-th row of $\F$),
we first note that $(\M^j)^*=\M^j$ and hence, in particular,
$M^j_{ik}=\overline{M^j_{ki}}$ (i.e., $\Re(M^j_{ik})=\Re(M^j_{ki})$
and $\Im(M^j_{ik})=-\Im(M^j_{ki})$). Thus, it is easily seen that for
each $i=1,\dots,N$, the terms in the double-sum $\x^\top
\M^j\x=\sum_{\l=1}^{N}\sum_{k=1}^{N} M^j_{\l k}x_{\l}x_{k}$ that
contain $\x_i$ are precisely
\[
M^j_{ii}x_i^2+x_i\sum_{k=1,k\neq i}^{N} M^j_{ik}x_k+x_i\sum_{\l=1,\l\neq i}^{N} M^j_{\l i}x_{\l}  = M^j_{ii}x_i^2 +\Big(2\sum_{k=1,k\neq i}^{N} \Re\big(M^j_{ik}\big)x_k\Big)x_i.
\]
Hence, $\tfrac{\partial}{\partial x_{i}}\big(y_j-\x^\top
\M^j\x\big)=-2\sum_{k=1}^{N}\Re(M^j_{ik})x_k=-2\Re(\M^j_{i\cdot})\x$,
and therefore $\nabla\phi(\x)$ is equal to
\[
-\frac{1}{2}\sum_{j=1}^N\big(y_j-\x^\top \M^j\x\big)\cdot\Big(2\Re(\M^j_{1\cdot})\x,\dots,2\Re(\M^j_{N\cdot})\x\Big)^\top = -\sum_{j=1}^{N}\big(y_j-\x^\top\F_{j\cdot}^*\F_{j\cdot}\x\big)\Re(\F_{j\cdot}^*\F_{j\cdot})\x.
\]
Consequently,
\begin{equation}\label{eq:altmodel_gradX}
  \nabla\xi(\x) = \nabla\phi(\x)+\nabla\psi(\x) =
  \mu\P_e^\top\big(\P_e\x-\hat\D\a\big)-\sum_{j=1}^{N}\big(y_j-\x^\top
  \F_{j\cdot}^*\F_{j\cdot}\x\big)\Re(\F_{j\cdot}^*\F_{j\cdot})\x.
\end{equation}
Thus, with the projection $\proj_{\cX}(\x)=\max\{0,\min\{1,\x\}\}$
(component-wise) and a step size $\gamma^{X}_{(\l)}>0$ (chosen by a
simple backtracking scheme to achieve a reduction in the
objective~\eqref{prob:PRDL1dXupdate}), the update of the phase
retrieval solution estimate in the $\l$-th DOLPHIn iteration for the 1D-case sets~$\x^{(l+1)}$ to the value
\[
\proj_{\cX}\Big(\x^{(\l)} - \gamma^{X}_{(\l)}\Big(\sum_{j=1}^{N}\big((\x^{(\l)})^\top \F_{j\cdot}^*\F_{j\cdot}\x^{(\l)}-y_j\big)\Re(\F_{j\cdot}^*\F_{j\cdot})\x^{(\l)} +\mu\P_e^\top\big(\P_e\x^{(\l)}-\hat\D^{(\l)}\a^{(\l)}\big)\Big)\Big).
\]

\begin{rem}
  The expression~\eqref{eq:altmodel_gradX} can be further simplified by rewriting $\nabla\phi(\x)$:
  Since the only complex-valued part within $\nabla\phi(\x)$ is $\F_{j\cdot}^*\F_{j\cdot}$, we can take
  the real part of the whole expression instead of just this product, i.e.,
  \[
    \nabla\phi(\x) 
    =-\sum_{j=1}^{N}\big(y_j-\x^\top\F_{j\cdot}^*\F_{j\cdot}\x\big)\Re(\F_{j\cdot}^*\F_{j\cdot})\x
    =
    \Re\Big(-\sum_{j=1}^{N}\big(y_j-\x^\top\F_{j\cdot}^*\F_{j\cdot}\x\big)\F_{j\cdot}^*\F_{j\cdot}\x\Big).
  \]
  Further, using $\F_{j\cdot}^*=(\F_{j\cdot})^*=(\F^*)_{\cdot j}$ and rearranging terms, this becomes
  \[
    \nabla\phi(\x)  =\Re\Big(\sum_{j=1}^{N}(\F^*)_{\cdot j}\big(\abs{\F\x}^2-\y\big)_j\F_{j\cdot}\x\Big)
     =\Re\Big(\F^*\big(\big(\abs{\F\x}^2-\y\big)\odot\F\x\big)\Big).
  \]
  From this last form, it is particularly easy to obtain the gradient
  matrix $\nabla\phi(\X)$ in the 2D-case, which corresponds precisely
  to the gradient $\nabla\phi(\x)$ with $\x$ the vectorized matrix
  variable $\X$ and $\F$ representing the linear operator $\cF(\X)$
  (in terms of $\x$), reshaped to match the size of $\X$ (i.e.,
  reversing the vectorization process afterwards). Similarly, the
  expression for $\nabla\psi(\X)$ can be derived from $\nabla\psi(\x)$
  w.r.t. the vectorized variable; here, the product with $\P_e^\top$
  then needs to be replaced by an application of the adjoint $\cE^*$,
  which can be recognized as $\cE^*(\cdot)=\mR\odot\cR(\cdot)$
  by translating the effect of multiplying by $\P_e^\top$ to the
  vectorized variable back to matrix form. (Analogously, one obtains
  $\cR^*(\mZ)=\cE((1/\mR)\odot\mZ)$, where $1/\mR$ has entries
  $1/r_{ij}$, i.e., is the entry-wise reciprocal of $\mR$.)
\end{rem}

\section{Convergence Proof}\label{sec:appendixB}
This section is devoted to proving our main Theorem~\ref{thm:convCGD}.
As mentioned in Section~\ref{subsec:monotone}, we can regard DOLPHIn
(Algorithm~\ref{alg:PRDL}) as a special case of the coordinate
gradient descent (CGD) method proposed in~\cite{TY09}. This method aims 
at solving problems of the form 
\begin{equation}\label{prob:cgd}
  \min_{\x} F_c(\x)\define f(\x)+cP(\x),
\end{equation}
where $c>0$, $P:\R^n\to(-\infty,+\infty]$ is proper, convex and lower
semicontinuous, and $f:\R^n\to\R$ is smooth on an open subset of
$\R^n$ containing $\text{dom\,}P=\{\x:P(\x)<\infty\}$. Assuming
(block-)separability of $P$, CGD then operates iteratively on
different blocks of variables (indexed by
$J_k\subseteq\{1,\dots,n\}$) by determining (in the
$k$-th iteration) a descent direction
\[
\z_k\define\arg\min_{\z}\Big\{\nabla f(\x_k)^\top \z+\tfrac12 \z^\top\H^k\z+cP(\x^k+\z)\,:\,{\z}_{{\bar{J}}_{k}}=\zeros\Big\},
\]
where $\H^k$ is such that $\H^k-\underline{\nu}\I$ and
$\bar\nu\I-\H^k$ are positive definite for all $\l$ (for some
constants $0<\underline{\nu}\leq\bar\nu$), and
$\bar{J}_k\define\{1,\dots,n\}\setminus J_k$, and carefully choosing a
step size $\alpha_k$ to set the new iterate
$\x_{\l+1}\define\x_k+\alpha_k\z_k$. In particular, for several
possible ways to define the sets $J_k$, (linear) convergence to a
stationary point can be guaranteed (under varying assumptions) when
the step sizes are chosen according to the following Armijo rule with
parameters $\beta,\sigma\in(0,1)$ and $\gamma\in[0,1)$: 

\begin{center}
\begin{minipage}{0.65\textwidth}
Choose
$\alpha_k^0>0$ and let $\alpha_k$ be the largest number in
$\{\alpha_k^0\beta^j\}_{j=0,1,2,\dots}$ such that
$F_c(\x_k+\alpha_k\z_k)-F_c(\x_k)\leq\alpha_k\sigma(\gamma-1)\z_k^\top\H^k\z_k$.
\end{minipage}
\end{center}

\noindent
For details, see Theorem~1, and (particularly) Sections~2, 5 and 8 in \cite{TY09}.

The following proof will explain how to relate DOLPHIn to the CGD
method and, consequently, establish analogous convergence results.
\begin{IEEEproof}[Proof of Theorem~\ref{thm:convCGD}]

  First, note that our problem~\eqref{prob:PRDL} can be cast in the
  form~\eqref{prob:cgd} via replacing the constraints $\X\in\cX$ and
  $\D\in\cD$ by adding the corresponding indicator functions
  $1_{\cX}(\X)$ and $1_{\cD}(\D)$ to the objective.  (For the sake of
  simplicity, we mostly omit explicitly converting between
  matrix-variables and vector-variables in this proof, and refer back
  to the corresponding explanations in the Appendix~\ref{sec:appendixA}.)

  We decompose the index set $J\define\{1,\dots,(p+s)n+N_1N_2\}$ for
  all components of $\A$, $\X$ and $\D$ into $p+n+1$ disjoint subsets:
  $J^{a,i}\define\{(i-1)n+1,\dots,in\}$ for the columns $\a^i$ of $\A$
  ($i=1,\dots,p$), $J^X\define\{pn+1,\dots,pn+N_1N_2\}$ for $\X$, and
  $J^{d,j}\define\{pn+N_1N_2+(j-1)s+1,\dots,pn+N_1N_2+js\}$ for the
  columns $\d^j$ of $\D$ ($j=1,\dots,n$). Each DOLPHIn iteration thus
  cycles through $J^{a,1},\dots,J^{a,p},J^X,J^{d,1},\dots,J^{d,n}$
  (in that order\footnote{The updates of $\a^i$ may again be performed
    in parallel, and could all use the same step sizes, so that one
    might alternatively consider $J^A\define\bigcup_{i=1}^{p}J^{a,i}$
    as one block; this amounts to Step~\ref{step:PRDL_Aupdate} as
    stated in Algorithm~\ref{alg:PRDL},
    cf.~Section~\ref{subsec:Aupdate}. For clarity, we stick to using $J^{a,i}$ here.})  and updates
  the associated blocks of variables. For convenience, we will
  henceforth index all occurring variables or parameters from the CGD
  perspective by the same notation as these index subsets, along with
  the DOLPHIn iteration counter $\l$; for instance, $\z^{a,i}_{\l}$
  will be the subvector (w.r.t. components indexed by $J^{a,i}$) of
  the CGD direction used for the update of $\a^i$, and $\H^{a,i}_{\l}$
  is the submatrix obtained from ``$\H^k$'' by restricting rows and
  columns to those indexed by~$J^{a,i}$. 

  Consider the $\l$-th DOLPHIn iteration. We start with the updates of
  $\a^i$ ($i=1,\dots,p$), which take the form of one ISTA iteration for
  \[
  \min_{\a}F_c^{a,i}(\a)\define\underbrace{\tfrac12\norm{\x^i-\D\a}_2^2}_{=:\xi^i(\a)}+\tfrac{\lambda}{\mu}\norm{\a}_1
  \]
  with $\x^i\equiv\x^i_{(\l)}$, $\D\equiv\D_{(\l)}$, and thus set
  \[
  \a^i_{(\l+1)}\define\cS_{\gamma^{a,i}_{\l}\lambda/\mu}\big(\a^i_{(\l)}-\gamma^{a,i}_{\l}\nabla\xi^i(\a^i_{(\l)})\big)=\cS_{\gamma^{a,i}_{\l}\lambda/\mu}\big(\a^i_{(\l)}-\gamma^{a,i}_{\l}\D^\top(\D\a^i_{(\l)}-\x^i)\big),
  \]
  with step size $\gamma^{a,i}_{\l}$ as described in
  Section~\ref{subsec:Aupdate}. From \cite[Lemma 2.3 (see also Remark
  3.1)]{BT09}, we know that
  \begin{equation}\label{eq:AupdArmijo}
  F_c^{a,i}(\a^i_{(\l+1)})-F_c^{a,i}(\a^i_{(\l)})\leq -\tfrac{1}{2\gamma^{a,i}_{\l}}\norm{\a^i_{(\l+1)}-\a^i_{(\l)}}_2^2.
  \end{equation}
  Note that the difference with respect to the original objective
  function of the DOLPHIn model~\eqref{prob:PRDL} actually is
  \[
  \mu F_c^{a,i}(\a^i_{(\l+1)})-\mu F_c^{a,i}(\a^i_{(\l)})=\mu(F_c^{a,i}(\a^i_{(\l+1)})-F_c^{a,i}(\a^i_{(\l)}))\leq F_c^{a,i}(\a^i_{(\l+1)})-F_c^{a,i}(\a^i_{(\l)}),
  \]
  since $F_c^{a,i}(\a^i_{(\l+1)})-F_c^{a,i}(\a^i_{(\l)})\leq 0$ and
  $\mu\geq 1$.  Therefore, the bound~\eqref{eq:AupdArmijo} also holds
  for the full combined DOLPHIn objective.

  Let $\sigma\in(0,1)$ and $\gamma\in[1-\tfrac{1}{2\sigma},1)$; so,
  $\sigma(\gamma-1)\geq-\tfrac12$.  Then, our update of $\a^i$
  corresponds to a CGD update w.r.t. the block $J^{a,i}$ with
  $\H^{a,i}_{\l}\define(\gamma^{a,i}_{\l})^{-1}\I_n$ and step size
  $(\alpha_k=)1$ that satisfies the Armijo rule (with arbitary
  $\alpha_k^0>0$, $\beta\in(0,1)$) stated at the beginning of this
  section \emph{if the CGD direction equals
    $\a^i_{(\l+1)}-\a^i_{(\l)}$}. Thus, we need to show that
  $\a^i_{(\l+1)}=\a^i_{(\l)}+\z^{a,i}_{\l}$ holds for
  \[
  \z^{a,i}_{\l}=\arg\min_{\z}\nabla\xi^i(\a^i_{(\l)})^\top\z+\tfrac12\z^\top\H^{a,i}_{\l}\z+\tfrac{\lambda}{\mu}\norm{\z+\a^i_{(\l)}}_1
  =\arg\min_{\z}\nabla\xi^i(\a^i_{(\l)})^\top\z+\tfrac{1}{2\gamma^{a,i}_{\l}}\norm{\z}_2^2+\tfrac{\lambda}{\mu}\norm{\z+\a^i_{(\l)}}_1.
  \]  
  Since $\H^{a,i}_{\l}$ is diagonal and the $\l_1$-norm separable, the
  $j$-th component of the solution is determined (see, e.g., \cite[p.~393]{TY09}) as
  \[
  -\midpt\big\{ \gamma^{a,i}_{\l}\big(\nabla\xi^i(\a^i_{(\l)})_j-\tfrac{\lambda}{\mu}\big), \big(\a^i_{(\l)}\big)_j, \gamma^{a,i}_{\l}\big(\nabla\xi^i(\a^i_{(\l)})_j+\tfrac{\lambda}{\mu}\big) \big\}
   =\begin{cases}-\big(\a^i_{(\l)}\big)_j, &\abs{\kappa^{a,i,\l}_j}<\tfrac{\gamma^{a,i}_{\l}\lambda}{\mu};\\
    -\gamma^{a,i}_{\l}\big(\nabla\xi^i(\a^i_{(\l)})_j-\tfrac{\lambda}{\mu}\big), &\kappa^{a,i,\l}_j<-\tfrac{\gamma^{a,i}_{\l}\lambda}{\mu};\\
    -\gamma^{a,i}_{\l}\big(\nabla\xi^i(\a^i_{(\l)})_j+\tfrac{\lambda}{\mu}\big), &\kappa^{a,i,\l}_j>\tfrac{\gamma^{a,i}_{\l}\lambda}{\mu},\end{cases}
  \]
  where $\kappa^{a,i,\l}_j\define\big(\a^i_{(\l)}-\gamma^{a,i}_{\l}\nabla\xi^i(\a^i_{(\l)})\big)_j$, i.e., we indeed obtain 
  \[
  \z^{a,i}_{\l} = \cS_{\gamma^{a,i}_{\l}\lambda/\mu}\big(\a^i_{(\l)}-\gamma^{a,i}_{\l}\nabla\xi^i(\a^i_{(\l)})\big)-\a^i_{(\l)} = \a^i_{(\l+1)}-\a^i_{(\l)}.
  \]
  (We emphasize again that the choice of $\gamma^{a,i}_{\l}$ can be
  made separately for each $\a^i$ or the same for all of them. In
  particular, since $\D_{(\l)}$ remains unchanged throughout the
  update of $\A$, we could use
  $\gamma^{a,i}_{\l}=\gamma^A_{\l}\in(0,1/\norm{\D_{(\l)}}_2^2)$ for
  all~$i$, and summarize the separate updates to the expression for
  the matrix-variable $\A_{(\l+1)}$ as in Step~\ref{step:PRDL_Aupdate}
  of Algorithm~\ref{alg:PRDL}.) 
  
  Turning to the $\X$-update, we can proceed quite similarly; for the
  sake of more easily relating to the CGD expressions, let us consider
  the update problem in terms of the vectorized variable
  $\x=\vect(\X)$ (reusing notation from the Appendix~\ref{sec:appendixA}):
  \[
  \min_{\x}F_c^x(\x)\define\underbrace{\tfrac14\norm{\y-\abs{\F\x}^2}_2^2+\tfrac{\mu}{2}\norm{\P_e\x-\hat\D\a}_2^2}_{=:f^{x}(\x)=\phi(\x)+\psi(\x)}+1_{\cX}(\x),
  \]
  with $\cX\equiv[0,1]^{N_1 N_2}$, $\hat\D\equiv\D_{(\l)}\otimes\I_p$
  and $\a\equiv\vect(\A_{(\l+1)})$. In DOLPHIn, we set
  \[
  \x_{(\l+1)}\define \proj_{\cX}\big(\x_{(\l)}-\gamma^X_{\l}\nabla f^x(\x_{(\l)})\big),
  \]
  with $\gamma^X_{\l}$ chosen by the Armijo rule stated in
  Theorem~\ref{thm:convCGD}, i.e., as the largest number in
  $\{\bar\eta\eta^j\}_{j=0,1,2,\dots}$ such that
  \[
  f^x\Big(\proj_{\cX}\big(\x_{(\l)}-\gamma^X_{\l}\nabla f^x(\x_{(\l)})\big)\Big)-f^x(\x_{(l)})
  \leq -\tfrac{1}{2\gamma^X_{\l}}\norm{\proj_{\cX}\big(\x_{(\l)}-\gamma^X_{\l}\nabla f^x(\x_{(\l)})\big)-\x_{(\l)}}_2^2.
  \]
  Similarly to the $\A$-update, this conforms with taking a step of
  length 1 that satisfies the CGD Arminjo rule mentioned earlier if
  the CGD descent direction associated to the $\x$-update, computed
  using $\H^x_{\l}\define(\gamma^X_{\l})^{-1}\I_{N_1N_2}$, is
  equal to $\x_{(\l+1)}-\x_{(\l)}$. To see that this is indeed the
  case, observe that
  \[
  \z^x_{\l}\define\arg\min_{\z}\nabla f^x(\x_{(\l)})^\top\z+\tfrac12\z^\top\H^x_{\l}\z+1_{\cX}(\z+\x_{(\l)})
  =\arg\min_\z\{\,\nabla f^x(\x_{(\l)})^\top\z+\tfrac{1}{2\gamma^X_{\l}}\norm{\z}_2^2\,:\,\z\in\cX\,\}
  \]  
  is again fully separable, and by analyzing the one-dimensional
  problems per component, it is easily verified that
  \[
  \z^x_{\l}=\proj_{\cX-\x_{(\l)}}\big(-\gamma^X_{\l}\nabla f^x(\x_{(\l)})\big)
  =\proj_{\cX}\big(\x_{(\l)}-\gamma^X_{\l}\nabla f^x(\x_{(\l)})\big)-\x_{(\l)}.
  \]
  
  (Note that, from the CGD viewpoint, the Armijo search in our
  $\X$-update is ``moved'' into the search direction computation; this
  variant and the possible increase in computational cost of the
  search is briefly commented on in \cite[Section 8]{TY09}. Moreover, here the 
  objective reduction w.r.t. $f^x$ is exactly that achieved for the original 
  DOLPHIn objective~\eqref{prob:PRDL}, since the relevant objective parts were 
  not rescaled.)

  Finally, consider the update of $\D$, which is done in a
  block-coordinate fashion in
  DOLPHIn. Steps~\ref{step:PRDL_DupdateStart}--\ref{step:PRDL_DupdateEnd}
  of Algorithm~\ref{alg:PRDL} solve
  \begin{equation}\label{prob:DupdateCol}
    \min_{\d}F_c^{d,j}(\d)\define\underbrace{\tfrac{1}{2}\sum_{i=1}^{p}\norm{\x^i-\sum_{\begin{subarray}{c}k=1\\k\neq j\end{subarray}}^{n}a^i_k\d^k-a^i_j\d}_2^2}_{=:f^{d,j}(\d)}+1_{\norm{\cdot}_2\leq 1}(\d)
  \end{equation}
  for each column $\d^j$ of $\D$ in turn (with
  $\x^i\equiv\x^i_{(\l+1)}$, $\a^i\equiv\a^i_{(\l+1)}$ and the
  remaining dictionary columns $\d^k\equiv(\d^k)_{(\l)}$, $k\neq j$).
  Indeed, these subproblems are solved \emph{exactly}, with the closed
  form solution easily derived as
  \[
  (\d^j)_{(\l+1)}\define \proj_{\norm{\cdot}_2\leq 1}\left(\frac{1}{\sum_{i=1}^{p}(a^i_j)^2}\sum_{i=1}^p a^i_j\big(\x^i-\sum_{\begin{subarray}{c}k=1\\k\neq j\end{subarray}}^{n}a^i_k\d^k\big)\right)
  =\frac{1}{\max\{1,\norm{\tfrac{1}{w_j}\q^j}_2\}}\big(\tfrac{1}{w_j}\q^j\big)
  \]
  with $w_j\define\sum_{i} (a^i_j)^2$ and $\q^j\define\sum_{i}
  a^i_j\big(\x^i-\sum_{k\neq j}a^i_k\d^k\big)$.  (If $w_j=0$, and thus
  $a^i_j=0$ for all $i=1,\dots,p$, then column $\d^j$ is not used in
  any of the current patch representations; in that case,
  problem~\eqref{prob:DupdateCol} has a constant objective and is
  therefore solved by any arbitrary $\d$ with $\norm{\d}_2\leq 1$,
  e.g., a normalized random vector as in
  Step~\ref{step:PRDL_DupdateReset} of Algorithm~\ref{alg:PRDL}.
  Since obviously, no descent can be achieved for a constant
  objective, the Armijo rule is not applicable in such cases, so in
  the remainder of the proof we consider only $j$ with $w_j\neq 0$.)

  Note that the update of $\d^j$ can also be seen as one ISTA iteration, 
  starting at $(\d^j)_{(\l)}$ and using a fixed step size of $1/w_j$; in fact, 
  $w_j$ is the smallest Lipschitz constant of the gradient $\nabla f^{d,j}(\d)=w_j\d-\q^j$:
  \[
  \norm{\nabla f^{d,j}(\d)-\nabla f^{d,j}(\hat\d)}_2=\norm{w_j(\d-\hat\d)}_2\leq w_j\norm{\d-\hat\d}_2
  \]
  by the Cauchy-Schwarz inequality. In ISTA terms, $(\d^j)_{(\l)}$ is
  computed as the (unique) minimizer of $1_{\norm{\cdot}_2\leq 1}(\d)$
  plus the quadratic approximation $f^{d,j}((\d^j)_{(\l)}+\nabla
  f^{d,j}((\d^j)_{(\l)})^\top(\d-(\d^j)_{(\l)})+\tfrac{1}{w_j}\norm{\d-(\d^j)_{(\l)}}_2^2$
  of $f^{d,j}(\d)$. Since $f^{d,j}$ is quadratic with Hessian
  $\nabla^2 f^{d,j}=w_j\I_s$, this approximation is precisely the
  Taylor expansion of $f^{d,j}(\d)$ around $(\d^j)_{(\l)}$, so the
  ISTA update coincides with ours. By this equivalence, we can
  immediately conclude from \cite[Lemma 2.3 (and Remark 3.1)]{BT09}
  that
  \[
     F_c^{d,j}((\d^j)_{(\l+1)})-F_c^{d,j}((\d^j)_{(\l)}) \leq-\tfrac{1}{2w_j}\norm{(\d^j)_{(\l+1)}-(\d^j)_{(\l)}}_2^2.
  \]
  Analogously to the $\A$-update, the term on the left hand side
  upper-bounds the corresponding difference w.r.t. the original
  DOLPHIn objective function~\eqref{prob:PRDL}, since $\mu\geq
  1$. Hence, we recognize again our desired Armijo rule, and therefore
  proceed to show that $\z^{d,j}_{\l}=(\d^j)_{(\l+1)}-(\d^j)_{(\l)}$
  when computed using $\H^{d,j}_{\l}=w_j\I_s$ (and CGD step size
  1). By definition,
  \[
  \z^{d,j}_{\l}\define\arg\min_{\z}\nabla f^{d,j}((\d^j)_{(\l)}^\top\z+\tfrac{w_j}{2}\z^\top\z+1_{\norm{\cdot}_2\leq 1}(\z+(\d^j)_{(\l)});
  \]
  substituting $\z=\d-(\d^j)_{(\l)}$ $\Leftrightarrow$
  $\d=\z+(\d^j)_{(\l)}$, we can alternatively solve
  \[
  \arg\min_{\d}\nabla f^{d,j}((\d^j)_{(\l)}^\top(\d-(\d^j)_{(\l)})+\tfrac{w_j}{2}\norm{\d-(\d^j)_{(\l)}}_2^2+1_{\norm{\cdot}_2\leq 1}(\d).
  \]
  Making use of the Taylor expansion of $f^{d,j}(\d)$ around $(\d^j)_{(\l)}$, this can be rewritten as
  \[
    \arg\min_{\d} f^{d,j}(\d)-\underbrace{f^{d,j}((\d^j)_{(\l)})}_{\text{constant}}+1_{\norm{\cdot}_2\leq 1}(\d)\quad
    {\Leftrightarrow}\quad\arg\min_{\d} f^{d,j}(\d)+ 1_{\norm{\cdot}_2\leq 1}(\d),
  \]
  which is precisely the $\d^j$-update
  problem~\ref{prob:DupdateCol}. Therefore,
  $\z^{d,j}_{\l}=(\d^j)_{(\l+1)}-(\d^j)_{(\l)}$ holds indeed.

  Having established these relations between DOLPHIn and CGD, the
  claim of Theorem~\ref{thm:convCGD} now follows directly from the
  ``relaxed Armijo search'' variant of \cite[Theorem~1(e)]{TY09} as
  discussed in Section~8 of that paper. (The boundedness requirements
  on the DOLPHIn step size parameters ensure the prerequisite Assumption~1 in
  \cite{TY09}.)
\end{IEEEproof}

It is worth mentioning that a linear rate of convergence can also be
proved for DOLPHIn with the Armijo rule from
Theorem~\ref{thm:convCGD}, by extending the results of
\cite[Theorem~2 (cf. Section~8)]{TY09}.

Moreover, the requirement $\mu\geq 1$ can be dropped if the update problems
  for $\A$ and $\D$ are \emph{not} modified (by dividing the relevant
  objective parts by $\mu$) to a slightly more convenient, or familiar, form. 
  Then, to update column $\a^i$, one ISTA step for
  \[
  \min_{\a}\tfrac{\mu}{2}\norm{\x^i-\D\a}_2^2+\lambda\norm{\a}_1
  \]
  amounts to setting
  \[
  \a^i_{(\l+1)}\define\cS_{\gamma^{a,i}_{\l}\lambda}\big(\a^i_{(\l)}-\gamma^{a,i}_{\l}\mu\D^\top(\D\a^i_{(\l)}-\x^i)\big),
  \]
  and the update problem for each respective $\d^j$,
  \[
  \min_{\d}\tfrac{\mu}{2}\sum_{i=1}^p\norm{\x^i-\sum_{\begin{subarray}{c}k=1,\\k\neq j\end{subarray}}a^i_k\d^k-a^i_j\d}_2^2+1_{\norm{\cdot}_2\leq 1}(\d),
  \]
  has a closed form solution that defines the next iterate, namely
  \[
  (\d^j)_{(\l+1)}\define \frac{1}{\max\{1,\norm{\tfrac{\mu}{w_j}\q^j}_2\}}\big(\tfrac{\mu}{w_j}\q^j\big),
  \]
  with $w_j$ and $\q^j$ as defined in the above proof. Then, the ISTA
  descent guarantees that lead to fulfilling the Armijo condition
  needed to apply the CGD convergence result directly apply to the
  original DOLPHIn objective function, cf.~\eqref{prob:PRDL} and (with
  minor changes regarding the directions), the above proof goes through 
  completely analogously.


\bibliographystyle{IEEEtran}
\bibliography{DOLPHIn}

\end{document}